\documentclass[10pt]{amsart}
\usepackage{rotating}
\usepackage{multirow}
\newcommand{\Z}{{\mathbb Z}}
\newcommand{\R}{{\mathbb R}}

\newcommand{\Q}{{\mathbb Q}}
\newcommand{\F}{{\mathbb F}}
\newcommand{\bbP}{{\mathbb P}}
\newcommand{\C}{{\mathbb C}}
\newcommand{\cK}{{\mathcal{K}}}

\newcommand{\Gal}{\mbox{Gal}}

\newcommand{\ord}{\mathop{\rm ord}}

\newcommand{\til}{\tilde} %widetilde other option

\newcommand{\fr}{\!\!\!\!\!\!/\!\!\!\!\!\!}
\newcommand{\cd}{\! \cdot \!}

\numberwithin{table}{section}
\numberwithin{equation}{section}

\title[Lightly ramified number fields with Galois group $S.M_{12}.A$]{
Lightly ramified number fields \\ with Galois group $S.M_{12}.A$}
\author{David P.\ Roberts}
\address{Division of Science and Mathematics, University of
  Minnesota Morris,  Morris, MN 56267}
\email{roberts@morris.umn.edu}
\begin{document}
\maketitle

\begin{abstract} We specialize various three-point covers to find number fields with Galois group 
$M_{12}$, $M_{12}.2$, $2.M_{12}$, or $2.M_{12}.2$ and light ramification in 
various senses.    One of our  $2.M_{12}.2$ fields has  the unusual property that it is ramified
only at the single prime $11$.  
\end{abstract}
      
 \section{Introduction}
 
    The Mathieu group $M_{12} \subset S_{12}$ is the second smallest of the twenty-six 
 sporadic finite simple groups, having order 
$\mbox{95,040} = 2^{6} \cdot 3^3 \cdot 5 \cdot 11$.  The outer automorphism group
 of $M_{12}$ has order $2$, and accordingly one has another interesting
 group $\mbox{Aut}(M_{12}) = M_{12}.2 \subset S_{24}$.  The Schur multiplier of
 $M_{12}$ also has order $2$, and one has a third interesting group
 $\tilde{M}_{12} = 2.M_{12} \subset S_{24}$.   Combining these last two extensions in the standard way,
 one gets a fourth interesting group $\tilde{M}_{12}.2 = 2.M_{12}.2 \subset S_{48}$.  
 
     In this paper we consider various three-point covers, some of which have appeared 
  in the literature previously \cite{MM,Mat,MZ}.  We
 specialize these three-point covers to get number fields with Galois group
 one of the four groups $S.M_{12}.A$ just discussed.  Some of these 
 number fields are unusually lightly ramified in various senses.    Of particular
 interest is a number field with Galois group $\tilde{M}_{12}.2$ ramified only at the
 single prime $11$.   Our general goal, captured by our title, it 
 to get as good a sense as currently possible of the most lightly
 ramified fields with Galois group $S.M_{12}.A$ as above.  
  
       Section~\ref{background} provides some general background information.
   Section~\ref{covers}
  introduces the three-point covers that we use.  Section~\ref{dessins} draws the
  dessins associated to these covers so that the presence of $M_{12}$ and
  its relation to $M_{12}.2$ can be seen very clearly.  Section~\ref{specialization}
  describes the specialization procedure.  
  Section~\ref{light} focuses on $M_{12}$ and $M_{12}.2$ 
  and presents number fields with small root discriminant, small Galois root discriminant, and small number of ramifying primes,
  These last three notions are related
  but inequivalent interpretations of ``lightly ramified.''  Finally
  Section~\ref{lifts} presents some explicit
   lifts to $\tilde{M}_{12}$ and  $\tilde{M}_{12}.2$.  
   
      We thank the Simons Foundation for partially supporting this work through grant \#209472.

\section{General background}
\label{background}

    This section provides general background information
to provide some context for the rest of this paper.

\subsection{Tabulating number fields}
 Let $G \subseteq S_n$ be a transitive permutation group of degree $n$, considered
up to conjugation.  Consider the set $\cK(G)$ of isomorphism classes of degree $n$ number fields $K$ with splitting field $K^g$ having Galois group $\Gal(K^g/\Q)$ equal to $G$.   The inverse
Galois problem is to prove that all $\cK(G)$ are non-empty.  The general expectation
is that all $\cK(G)$ are infinite, except for the special case $\cK(\{e\}) = \{\Q\}$.  

To study fields $K$ in $\cK(G)$, it is natural to focus on their discriminants
$d(K) \in \Z$.     A fundamental reason to focus on discriminants is that
the prime factorization $\prod p^{e_p}$ of $|d(K)|$ measures 
by $e_p$ how much any given prime $p$ ramifies 
in $K$.   In a less refined way, the size $|d(K)|$ is a measure of the complexity of $K$.  
In this latter context, 
to keep numbers small and facilitate comparison
between one group and another, it is generally better to work with 
the root discriminant $\delta(K) = |d(K)|^{1/n}$.  
     
     To study a given $\cK(G)$ computationally, a methodical approach is to 
explicitly identify the subset $\cK(G,C)$ consisting of all fields with root
discriminant at most $C$ for as large a cutoff $C$ as possible.   Often
one restricts attentions to classes of fields which are of particular interest, 
for example fields with $|d(K)|$ a prime
power, or with $|d(K)|$ divisible only by a prescribed
set of small primes, or with complex conjugation
sitting in a prescribed conjugacy class $c$ of $G$. 
All three of these last
conditions depend only on $G$ as an abstract group, not
on the given permutation representation of $G$.   In this
spirit, it is natural to focus on the Galois root discriminant $\Delta$
of $K$, meaning the root discriminant of $K^g$.  
One has $\delta \leq \Delta$. To fully compute
 $\Delta$, one needs to identify the inertia subgroups $I_p \subseteq G$
 and their filtration by higher ramification groups. 

     Online tables associated to \cite{JRNF} and \cite{KM} provide a large
 amount of information on low degree number fields.  The tables for
 \cite{JRNF} focus
 on completeness results in all the above settings,
 with almost all currently posted 
 completeness results being in
 degrees $n \leq 11$.    The tables for \cite{KM} cover
 many more groups as they contain 
 at least one field for almost every pair $(G,c)$ in
 degrees $n \leq 19$. For
 each $(G,c)$, the field with the smallest
 known $\delta$ is highlighted.  

  There is an increasing sequence
of numbers $C_1(n)$ such that $\cK(G,C_1(n))$ is known to be
empty by discriminant bounds for all $G \subseteq S_n$.  
Similarly, if one assumes
the generalized Riemann hypothesis, there are larger numbers
$C_2(n)$ for which one knows $\cK(G,C_2(n))$ is non-empty.  In the limit of
large $n$, these numbers tend to $4 \pi e^\gamma \approx 22.3816$
and $8 \pi e^\gamma \approx 44.7632$ respectively.    This last constant
especially is  useful as a reference point when considering
root discriminants and Galois root discriminants.   See e.g.\ \cite{Mar} for 
explicit instances of these numbers $C_1(n)$ and $C_2(n)$.

Via class field theory, identifying $\cK(G,C)$ for any
solvable $G$ and any cutoff $C$ can be regarded as a 
computational problem.  For $G$ abelian, one
has an explicit description of $\cK(G)$ in its entirety.  
For many non-abelian solvable $G$ one 
can completely identify very large $\cK(G,C)$. 
Identifying $\cK(G,C)$ for nonsolvable groups
is also in principle a computational problem.  
However run times are prohibitive in general and
only for a very limited class
of groups $G$ have non-empty $\cK(G,C)$
been identified.  

\subsection{Pursuing number fields for larger groups.} 
\label{pursuing}  When producing complete non-empty lists for a given
 $G$ is currently infeasible, one would nonetheless
like to produce as many lightly ramified fields as possible.  One can
view this as a search for best fields in $\cK(G)$ in various senses.
Here our focus is on smallest root discriminant $\delta$, smallest Galois root discriminant $\Delta$, and
smallest $p$ among fields
ramifying at a single prime $p$.

For the twenty-six sporadic groups $G$ in their smallest permutation 
representations, the situation is as follows.   The set $\cK(G)$ is known to be
infinite for all groups except for the Mathieu group 
$M_{23}$, where it is not even known to be non-empty \cite{MM}.  
One knows very little about ramification in these fields.   
Explicit polynomials are known only for $M_{11}$, $M_{12}$, 
$M_{22}$, and $M_{24}$.  The next smallest degrees 
come from the Hall-Janko group $HJ$ and Higman-Sims  
group $HS$, both in $S_{100}$.  The remaining
sporadic groups seem well beyond current reach in
terms of explicit polynomials because of their large degrees. 

For $M_{11}$, $M_{12}$, $M_{22}$, and $M_{24}$, one knows
infinitely many number fields, by specialization from a small
number of parametrized families.   In terms of known
lightly ramified fields, the situation is different for
each of these four groups.
The known $M_{11}$ fields come from
specializations of 
$M_{12}$ families satisfying certain strong conditions 
 and so instances with small discriminant
are relatively rare.   On \cite{KM}, the current records for smallest root discriminant
are give by the polynomials
\begin{eqnarray*}
f_{11}(x) & = & x^{11}+2 x^{10}-5 x^9+50 x^8+70 x^7-232 x^6+796 x^5+1400 x^4 \\ && 
 \qquad -5075 x^3+
10950 x^2+2805 x-90,  \\
&&\\
f_{12}(x) & = & x^{12}-12 x^{10}+8 x^9+21 x^8-36 x^7+192 x^6-240 x^5-84 x^4 \\ && \qquad 
+68 x^3-72
    x^2+48 x+5.
\end{eqnarray*}
The respective root discriminants are
\[
\begin{array}{rclcl}
\delta_{11} & = & (2^{18} 3^8 5^{18})^{1/11} &\approx & 96.2, \\
\delta_{12} & = & (2^{24} 3^{12} 29^{4})^{1/12} &  \approx & 36.9.
\end{array}
\]
Galois root discriminants
are much harder to compute in general, with the general method being sketched in 
\cite{JRGRD}.   The interactive website \cite{JRLF} greatly facilitates GRD computations, 
 as indeed in favorable cases it  computes
GRDs automatically.    In the two current cases the
GRDs are respectively
\[
\begin{array}{rclcl}
\Delta_{11} & = & 2^{13/6} 3^{7/8} 5^{39/20} & \approx & 270.8 \\ 
 \Delta_{12} & = &  2^{43/16} 3^{25/18} 29^{1/2} & \approx & 159.4. 
\end{array}
\]
  The ratios $96.2/36.8 \approx 2.6$ and 
$270.8/159.5 \approx 1.7$ are large already, 
especially considering the fact that $M_{12}$ is twelve times as large as $M_{11}$.
But, moreover, the sequence of known root discriminants increases much
more rapidly for $M_{11}$ than it does for $M_{12}$.  There is 
one known family each for $M_{22}$ \cite{Mal22} and $M_{24}$ \cite{Gr, Mu}.  
The $M_{22}$ family gives some specializations with root discriminant of order
of magnitude similar to those above.  The $M_{24}$ family seems to 
give fields only of considerably larger root discriminant.  

In this paper we focus 
not especially on $M_{12}$ itself, but more so on its extension $M_{12}.2$, for
which more good families are available.     On the one
hand, we go much further than one can at present for any other extension $G.A$ of 
a sporadic simple group.   On the other hand, we expect that there are
many $M_{12}$ and $M_{12}.2$ fields of comparably light ramification
that are not accessible by our approach.

\subsection{$M_{12}$ and related groups} 
\begin{table}[htb]
\[
 \begin{array}{| lrl| ll | ll | rr |}
 \hline
          C & |C| & \mbox{freq} & \lambda_{12} & \lambda^t_{12} & \lambda_{24} & \lambda^t_{24} & \mbox{\S\ref{lifts12}} &
           \mbox{\S\ref{lifts24}} \\
 \hline
                   1A & 1 & 1/190080 & 1^{12} & 1^{12} & 1^{24} & 1^{24} & 1 & 0\\
                   1A & 1 & 1/190080 & 1^{12} & 1^{12} & 2^{12} & 2^{12} & 0 &1 \\
                   2A & 792 & 1/240 & 2^6 & 2^6 & 4^6 & 4^6  & 768 & 789  \\
                   2B & 495 & 1/384 & 2^4 1^4 & 2^4 1^4 & 2^{12} & 2^{12} & 470 & 503    \\
                   2B & 495 & 1/384 & 2^4 1^4 & 2^4 1^4 & 2^{8} 1^{8} & 2^{8} 1^8 &521 & 515    \\
                   3A & 1760 & 1/108 & 3^3 1^3 & 3^3 1^3 & 3^6 1^6 & 3^6 1^6 & 1735 & 1776 \\
                   3A & 1760 & 1/108 & 3^3 1^3 & 3^3 1^3 & 6^3 2^3 & 6^3 2^3 & 1823 & 1781 \\
                   3B & 2640 & 1/72 & 3^4 & 3^4 & 3^8 & 3^8 & 2702 & 2578 \\
                   3B & 2640 & 1/72 & 3^4 & 3^4 & 6^4 & 6^4 & 2649 &  2510\\
                   4A & 5940 & 1/32 & 4^22^2 & 4^21^4 & 4^4 2^4 & 4^4 2^2 1^4 & 6002 & \multirow{2}{*}{11992}  \\
                   4B & 5940 & 1/32 & 4^21^4 & 4^22^2 & 4^4 2^2 1^4 & 4^4 2^4 & 5993 &  \\
                   5A & 9504 & 1/20 & 5^21^2 & 5^21^2 & 5^4 1^4 & 5^4 1^4 & 9329 & 9415 \\
                   5A & 9504 & 1/20 & 5^21^2 & 5^21^2 & 10^2 2^2 & 10^2 2^2 & 9405 & 9613 \\
                   6A & 15840 & 1/12 & 6^2 & 6^2 & 12^2 & 12^2 & 15798 &  15819\\
                   6B & 15840 & 1/12 & 6321 & 6321 & 6^2 3^2 2^2 1^2 & 6^2 3^2 2^2 1^2  & 15863 &15590\\
                   6B & 15840 & 1/12 & 6321 & 6321 & 6^3 2^3 & 6^3  2^3  & 15881 & 15828 \\
                   8A & 23760 & 1/8 & 84 & 821^2 & 8^2 4^2 & 8^2 4 \; 2 1^2 & 23613 & \multirow{2}{*}{47707}\\
                   8B & 23760 & 1/8 & 821^2 & 84 & 8^2 4 \; 2 1^2 & 8^2 4^2 & 24022 &   \\
                   10A & 19008 & 1/10 & (10)2 & 102 & (20)4 & (20)4 & 19048 & 18965 \\
                   11AB & 17280 & 1/11 & (11)1 & (11)1 & 11^2 1^2 & 11^2 1^2 & 17031 & 17308 \\
                   11AB & 17280 & 1/11 & (11)1 & (11)1 & (22)2 & (22)2 & 17425 & 17194  \\
                   \hline
                   2C & 1584 &  1/120 & \multicolumn{2}{c|}{2^{12}} & \multicolumn{2}{c|}{2^{24}} &  & 1650 \\
                   4C & 7920 & 1/24 &   \multicolumn{2}{c|}{4^4 2^4 } & \multicolumn{2}{c|}{ 4^8 2^8 } &  & 7964 \\
                   4D & 15840 & 1/12 &   \multicolumn{2}{c|}{ 4^6} & \multicolumn{2}{c|}{ 8^6 } &  &15688  \\
                   6C & 31680 & 1/6 &  \multicolumn{2}{c|}{6^4 } & \multicolumn{2}{c|}{6^8  } &  &  31651\\
                   10BC & 38016 & 1/5 &   \multicolumn{2}{c|}{10^2 2^2 } & \multicolumn{2}{c|}{10^4 2^4  } &  &38245  \\
                   12A & 31680 & 1/6 &    \multicolumn{2}{c|}{ 12^2} & \multicolumn{2}{c|}{ 24^2  } &  & 31577 \\
                   12BC & 63360 & 1/3 &   \multicolumn{2}{c|}{ 12 \; 6 \; 4 \; 2 } & \multicolumn{2}{c|}{ 12^2 6^2 4^2 2^2 } &  & 63493 \\ 
                   \hline
                  \end{array}
 \]
 \caption{\label{conjclasses} 
 First seven columns: information on conjugacy classes of $S.M_{12}.A$ and their sizes.  Last two columns:
 distribution
 of factorization partitions $(\lambda_{12},\lambda_{12}^t,\lambda_{24},\lambda_{24}^t)$ 
 of polynomials $(f_{B},f_{B^t},\tilde{f}_{B},\tilde{f}_{B^t})$ from \S\ref{lifts12};
 distribution of factorization partitions $(\lambda_{12}\lambda_{12}^t,\lambda_{24}\lambda_{24}^t)$  of polynomials
 $(f_{D2},\tilde{f}_{D2})$ from \S\ref{lifts24}}
 \end{table}

To carry out our exploration, we freely use group-theoretical facts
about $M_{12}$ and its extensions.   Generators of $M_{12}$ and 
$M_{12}.2$ are given pictorially in Section~\ref{dessins}
and lifts of these generators to $\til{M}_{12}$ and
$\til{M}_{12}.2$ are discussed in Section~\ref{lifts}.  The Atlas \cite{Atlas} 
as always provides a concise reference for group-theoretic facts.  
Several sections of  \cite{CS} provide further useful background information,
making the very beautiful nature of $M_{12}$ clear.  To get 
a first sense of $M_{12}$ and its extensions, a understanding of conjugacy 
classes and their sizes is particularly useful, and information is given in 
Table~\ref{conjclasses}.   

To assist in reading Table~\ref{conjclasses}, note 
that $M_{12}$ has fifteen conjugacy classes, $1A$, \dots, $10A$, $11A$, $11B$
in Atlas notation.  All classes are rational except 
for $11A$ and $11B$ which are conjugate over 
$\Q(\sqrt{-11})$.    Three pairs of these classes become
one class in $M_{12}.2$, the new merged classes being 
$4AB$, $8AB$, and $11AB$.  Also there are nine entirely 
new classes in $M_{12}.2$, all rational except for 
the Galois orbits $\{10B,10C\}$ and $\{12B,12C\}$.  
The cover $\tilde{M}_{12}$ has $26$ conjugacy classes,
with all classes rational except for the
Galois orbits $\{8A1',8A1'' \}$, $\{8B1',8B1'' \}$, $\{10A2',10A2''\}$, $\{11A1,11B1\}$,
$\{11A2,11B2\}$.   The $21$ Galois orbits correspond
to the $21$ lines of Table~\ref{conjclasses} above the 
dividing line.   Finally $\tilde{M}_{12}.2$ has $34$ 
conjugacy classes, $20$ coming from the $26$ conjugacy
classes of $\tilde{M}_{12}$ and fourteen new ones.  
The seven lines below the divider give 
a quotient set of these new fourteen classes;
the lines respectively correspond to
$1$, $2$, $2$, $1$, $2$, $2$, and $4$ classes.

The columns $\lambda_{12}$, $\lambda_{12}^t$, $\lambda_{24}$,
$\lambda_{24}^t$ contain partitions of 
$12$, $12$, $24$, and $24$ respectively.    It is 
through these partitions that we see conjugacy classes
in $S.M_{12}.A$, either via cycle partitions of 
permutations or degree partitions of 
factorizations of polynomials into irreducibles
in $\F_p[x]$.    For $M_{12}$, the partition 
$\lambda_{12}$ corresponds to the given
permutation representation while $\lambda_{12}^t$ 
corresponds to the twin degree twelve permutation
as explained in \S\ref{G2} below.  Similarly for $\til{M}_{12}$,
one has $\lambda_{24}$ corresponding
to the given permutation representation
and $\lambda_{24}^t$ corresponding to 
its twin.   For $M_{12}.2$ one has only
the partition $\lambda_{12} \lambda_{12}^t$ of $24$,
For $\til{M}_{24}.2$, one likewise
has only the partition $\lambda_{24} \lambda_{24}^t$
of $48$.  

The existence of the biextension $2.M_{12}.2$ is part of the
exceptional nature of $M_{12}$.   In fact, the
 outer automorphism group of $M_n$  has order $2$ for 
$n \in \{12,22\}$ and has order $1$ for the other possibililities,
$n \in \{11,23,24\}$.    Similarly,
the Schur multiplier of $M_{12}$ and $M_{22}$ has
order $2$ and $12$ respectively, and order $1$ for 
the other $M_n$. The rest of this section consists of general comments, 
illustrated by contrasting $2.M_{12}.2$ with
$2.M_{22}.2$.    

\subsection{$G$ compared with $G.2$} 
\label{G2} 
For a group $G \subseteq S_n$ 
and a larger group $G.2$, there are two possibilities: either
the inclusion can be extened to $G.2$ or it can not.  
In the latter case, certainly $G.2$ 
embeds in $S_{2n}$, although it might
also embed in a smaller $S_m$, as is the
case for e.g.\ $S_6.2 \subset S_{10}$. 

The extension $M_{22}.2$ embeds in $S_{22}$ while
$M_{12}.2$ only first embeds in $S_{24}$.   
The fact that $M_{12}.2$ does not have a smaller
permutation representation perhaps is 
a reason for its relative lack of presence in
the explicit literature on the inverse Galois problem.

When $G.2$ does not embed in $S_n$, there is an
associated twinning phenomenon:  fields in
$\cK(G)$ come in twin pairs, with two twins $K_1$
and $K_2$ sharing a common splitting field 
$K^g$.    When the outer automorphism group
acts non-trivially on the set of conjugacy 
classes of $G$, then twin fields $K_1$ and
$K_2$ do not necessarily have to have the 
same discriminant; this is the case for 
$M_{12}$.  

\subsection{$G$ compared with $\tilde{G}$}  For 
a group $G \subseteq S_n$ and a double cover 
$\tilde{G}$, finding the smallest $N$ for which $\tilde{G}$ embeds in
$S_N$ can require an exhaustive analysis of subgroups.  
One needs to find a 
subgroup $H$ of $G$ of smallest possible index $N$ 
which splits in $\tilde{G}$ in the sense that
there is a group $\hat{H}$ in $\tilde{G}$ 
which maps bijectively to $H$.  The
subgroup $H$ also needs to satisfy the condition
that its intersection with all its conjugates in $G$ is trivial.    

In the case of $G=A_n$ and $\tilde{G} = \tilde{A}_n$, the
Schur double cover, the desired $N$ is typically much larger
than $n$.   Similarly $\tilde{M}_{22}$ first embeds
in $S_{352}$ and $\tilde{M}_{22}.2$ first embeds
in $S_{660}$.   The fact that one has the low degree
embeddings $\tilde{M}_{12} \subset S_{24}$ and
$\tilde{M}_{12}.2 \subset S_{48}$ greatly 
facilitates the study of $\cK(\tilde{M}_{12})$ and $\cK(\tilde{M}_{12}.2)$
via explicit polynomials.   These embeddings
arise from the fact that $M_{11}$ 
splits in $\tilde{M}_{12}$.  

\subsection{The nonstandard double extension $(2.M_{12}.2)^*$} 
\label{isocline} There
is a second non-split double cover $(2.M_{12}.2)^*$ of  $M_{12}.2$.  
We refer to the group $2.M_{12}.2$ we are working with throughout this paper as
the standard double cover, since the ATLAS \cite{Atlas} prints its
character table.  The isoclinic \cite[\S6.7]{Atlas} variant $(2.M_{12}.2)^*$ is
considered briefly in \cite{Br}, where the embedding 
$(2.M_{12}.2)^* \subset S_{48}$ is also discussed.  

   Elements of $2.M_{12}.2$ above elements in the
class $2C \subset M_{12}.2$ have cycle type $2^{24}$,
as stated on Table~\ref{conjclasses}.  In contrast,
elements in $(2.M_{12}.2)^*$ above elements in $2C$ have
cycle type $4^{12}$.   This different behavior plays 
an important role in \S\ref{lack}.

\section{Three-point covers}
\label{covers}

\subsection{Six partition triples}
     Suppose given three conjugacy classes $C_0$, $C_1$, $C_\infty$ in a centerless group $M \subseteq S_n$.
Suppose each of the $C_t$ is {\em rational} in the sense that whenever $g \in C_t$ and
$g^k$ has the same order as $g$, then $g^k \in C_t$ too.   Suppose that the triple
$(C_0,C_1,C_\infty)$ is {\em rigid} in the sense that there exists a unique-up-to-simultaneous-conjugation triple $(g_0,g_1,g_\infty)$ with $g_t \in C_t$, $g_0 g_1 g_\infty = e$, and $\langle g_0,g_1,g_\infty \rangle = M$.   Then the theory of three-point covers applies in its
simplest form: there exists a canonically defined cover degree $n$ cover $X$ of $\bbP^1$, ramified
only above the three points $0$, $1$, and $\infty$, with local monodromy
class $C_t$ about $t \in \{0,1,\infty\}$ and global monodromy $M$.  Moreover,
this cover is defined over $\Q$ and the set $S$ at which it has bad reduction
satisfies 
\begin{equation*}
S_{\rm loc} \subseteq S \subseteq S_{\rm glob}.
\end{equation*}
Here $S_{\rm loc}$ is the set of primes dividing the order of one of the elements in
a $C_t$, while $S_{\rm glob}$ is the set of primes dividing $|M|$.   

    An interesting fact about the $M_n$ is that they contain no rational rigid triples 
 $(C_0,C_1,C_\infty)$.   Accordingly, we will not be using the theory of
 three-point covers in its very simplest form.  Instead, for each of our $M_{12}$ covers there
 is a complication, always involving the number $2$, but in different 
 ways.  We will not be formal about how the general theory needs to be modified, as our
 computations are standard, and all we need is the explicit equations that 
 we display  below to proceed with our construction of number fields.
 
 We use the language of partition triples rather than class triples. 
 The only essential difference is that the two conjugacy classes $11A$ and 
 $11B$ give rise to the same partition of twelve, namely $(11)1$.  The
 six partition triples we use are listed in Table~\ref{sixtriples}.   As we will 
 see by direct computation, the sets $S$ of bad reduction are always of the form 
 $\{2,3,q\}$, thus strictly smaller than $S_{\rm glob} = \{2,3,5,11\}$.  
The extra prime $q$ is $5$ for Covers~$A$, $B$, and $B^t$, while
it is $11$ for Covers $C$, $D$, and $E$.

\begin{table}[htb]
\[
{\renewcommand{\arraycolsep}{4pt}
   \begin{array}{| llll | cccc |cccc|}
   \hline
   \mbox{Name} & \lambda_0 & \lambda_1 & \lambda_\infty & M_{12} & 
    M_{12}.2 &   \tilde{M}_{12} & \tilde{M}_{12}.2 & 2  &  3 & 5  & 11  \\
   \hline
    A & 3333 & 22221111 & (10)2 & &   \surd && &W & U & T &     \\
    B & 441111 & 441111 & (10)2 & \surd && \surd &  &  U & U & T &  \\
    B^t & 4422 & 4422 & (10)2 & \surd && \surd &  &  U & U & T &  \\
    C & 333111 & 222222 & (11)1 & &  \surd && \surd  & U & U & & T \\
    D& 3333 & 22221111 & (11)1 & &   \surd && \surd  & U & U & & T \\
    E & 333111 & 333111 & 66 & \surd &   \surd && &  W & T & & U    \\
    \hline
   \end{array}
}
\]
\caption{\label{sixtriples} Left: The six dodecic partition triples pursued
in this paper.  Middle:  The Galois groups $G$ they give rise to. 
Right: The primes of bad reduction 
and their least ramified behavior (Unramified, Tame, Wild) for specializations,
according
to Tables~\ref{abtable} and \ref{cdetable}.}
\end{table}

\subsection{Cover $A$} 
\label{CoverA} Cover~$A$
 was studied by Matzat \cite{Mat} in one of the first computational 
successes of the theory of three-point covers.  
The complication here is that there are two conjugacy
classes of $(g_0,g_1,g_\infty)$.   It turns out that they
are conjugate to each other over $\Q(\sqrt{-5})$.   Abbreviating
$a = \sqrt{-5}$, one finds  an equation for this
cover to be
\begin{eqnarray*}
f_{A}(t,x) & = &  5^3 \left(24 a x^2+16 a x-648 a-x^4-60
    x^3-870 x^2-220 x+6399\right)^3 \\
    & & -2^{12} 3^{15} (118 a-475) t x^2.
    \end{eqnarray*}
To remove irrationalities, we define
\begin{eqnarray*}
f_{A2}(t,x) &= & f_{A}(t,x) \overline{f}_A(t,x), 
\end{eqnarray*}
where $\overline{\cdot}$ indicates conjugation on coefficients.   
Because all cuspidal partitions
involved are stable under twinning, the generic Galois group
of $f_{A2}(t,x)$ is $M_{12}.2$, not the 
$M_{12}^2.2$ one might expect from similar situations in 
which quadratic irrationalities are removed in 
the same fashion.  

\subsection{Covers $B$ and $B^t$}
\label{CovsB} Cover $B$ is the most 
well-known of the covers in this paper, having
been introduced by Matzat and Zeh-Marschke \cite{MZ} and studied 
further in the context of lifting by Bayer, Llorente, and Vila \cite{BLV}
and Mestre \cite{Me}.    The complication from Cover~$A$ of there being 
two classes of $(g_0,g_1,g_\infty)$ is present here too. 
However, in this case, the complication can be addressed
without introducing irrationalities.  Instead one uses $\lambda_0=\lambda_1$
and twists accordingly.   An equation is then
\begin{eqnarray*}
f_{B}(s,x) & = & 3 x^{12}+100 x^{11}+1350 x^{10}+9300 x^9+32925 x^8+45000 x^7-43500
    x^6 \\ && -147000 x^5+46125 x^4+172500 x^3-16250 x^2+22500 x+1875 \\ 
    && - s 2^{11} 5^2 x^2.
\end{eqnarray*}
The twisting is seen in the polynomial discriminant, which is 
\[
D_{B}(s) = 2^{144} 3^{120} 5^{38} (s^2-5).  
\]
 So here and for $B^t$ below, the three critical values of the
cover are $-\sqrt{5}$, $\sqrt{5}$, and $\infty$.  The three critical values 
of Covers $A$, $C$, $D$, and $E$ are all at their standard positions 
$0$, $1$, and $\infty$.

    While the outer automorphism of $M_{12}$ fixes the conjugacy class $10A = (10)2$, it
switches the classes $4B=441111$ and $4A=4422$.  Therefore Cover $B^t$, 
the twin of 
Cover $B$, has ramification triple $(4422,4422,(10)2)$ and hence genus two.  While the
other five covers have genus zero and were easy to compute directly, it would be 
difficult to compute Cover $B^t$ directly.  Instead we started from $B$ and applied 
resolvent constructions, eventually ending at the following polynomial:
\begin{eqnarray*}
f_{B^t}(s,x) & = &  5^2 \left(2500 x^{12}-45000 x^{10}+310500 x^8-1001700 x^6+1433700
    x^4 \right. \\ && \qquad \left. -641520 x^2+174960 x+88209\right) \\
    && -270 s (12 x+25) \left(50 x^6-450 x^4+1080 x^2-297\right) \\
    && +3^6 s^2 (12 x-25)^2.
\end{eqnarray*}
Unlike in our equations for Covers $A$, $B$, $C$, and $D$, here 
$x$ is not a coordinate on the covering curve.  Instead the
covering curve is a desingularization of the plane curve
given by $f_{B^t}(s,x)=0$.  The function $x$ has degree two and there
is a degree $5$ function $y$ so that the curve can be presented
in the more standard form $y^2 = 15 x \left(5 x^4+30 x^3+51 x^2-45\right)$.

\subsection{Covers $C$ and $D$} 
\label{CovsCD} The next two covers are remarkably similar to each other 
and we treat them simultaneously.  In these cases, the underlying permutation 
triple is rigid.  However it is not rational, since $11A$ and $11B$ are conjugates are
one another over $\Q(\sqrt{-11})$.  So as for Cover $A$, there is an irrationality in our 
final polynomials, although this time we knew before computing that the field of
definition would be $\Q(\sqrt{-11})$.  Abbreviating $u = \sqrt{-11}$, our polynomials are
\begin{eqnarray*}
f_{C}(t,x) & = &  \left(21 u
    x+13 u-2 x^3-54 x^2-321 x-83\right)^3 \\
&&     \cdot  \left(69 u x+1573 u-2 x^3-102 x^2-1713 x-10043\right) \\
 &&   + t2^9 3^{12} (253 u+67) t x, \\
f_{D}(t,x) & = & -11^2 u \left(1188 u x^3+198 u x^2-1346 u x-27 u+594 x^4 \right. \\
&& \qquad \left. -7920
    x^2-1474 x+135\right)^3 \\
    && - 2^8 3^{13} (253 u-67) t x.
 \end{eqnarray*}
 As with Cover $A$, we remove irrationalities by forming 
 $f_{C2}(t,x) =  f_{C}(t,x) \overline{f}_C(t,x)$ and  
 $f_{D2}(t,x) =  f_{D}(t,x) \overline{f}_D(t,x)$.  As for 
 $f_{A2}(t,x)$, the Galois group of these new polynomials
 is $M_{12}.2$.

\subsection{Cover $E$}  
\label{CovE} 
One complication for Cover $E$ is the same as for Covers $A$, $B$, and $B^t$:  there are two classes of underlying
$(g_0,g_1,g_\infty)$.    As for $B$ and $B^t$, the classes $C_0$ and
$C_1$ agree, which can be exploited by twisting to obtain rationality.   But
now, unlike for $B$ and $B^t$, this class, namely $3A$, is
stable under twinning.  So now, replacing the twin pair $(X_B, X_{B^t})$, there is a single curve 
$X_E$ with a self-twinning involution.
Like $X_B$, this curve  has
genus zero and is defined over $\Q$.  However, 
a substantial complication
arises only here: the curve $X_E$ does not have a rational point  
and is hence not parametrizable.

    We have computed a corresponding degree twelve polynomial $f_E(s,x)$ and used it to determine
 a degree twenty-four polynomial  
   \begin{eqnarray*}
   \lefteqn{f_{E2}(t,x) =} \\
&& (1-t) \left(x^6-20 x^5+262 x^4-15286 x^3+477665 x^2-10170814
    x+96944940\right)^3 \\ 
    &&  \left(x^6+60 x^5+2406 x^4+56114 x^3+1941921
    x^2+55625130 x+996578748\right) \\
    &&+t \left(x^{12}+396 x^{10}-27192 x^9+933174
    x^8-20101752 x^7+169737744 x^6 \right.  \\ && \qquad  -16330240872 x^5+538400028969
    x^4-8234002812376 x^3 \\ && \qquad \left. +195276967064388 x^2-3991355037576144
    x+30911476378259268\right)^2 \\
    &&+2^4 3^{12} 11^{22} t (t-1).
    \end{eqnarray*}
  One recovers $f_E(s,x)$ via $f_{E2}(1+s^2/11,x) = f_E(s,x) f_E(-s,x)$.  The 
  discriminants of $f_E(s,x)$ and $f_{E2}(t,x)$ are respectively 
  \begin{eqnarray*}
  D_E(s) & = & 2^{64} 3^{48} 11^{60} (s^2+11)^6 c_4(s)^2, \\
  D_{E2}(t) & = & 2^{224} 3^{168} 11^{264} t^{12}  (t-1)^{12}  c_{10}(t).^2
  \end{eqnarray*}
  The last factors in each case have the indicated degree and do not contribute to 
  field discriminants.  The Galois group of $f_{E}(s,x)$ over $\Q(s)$ is $M_{12}$
  and $f_E(-s,x)$ gives the twin $M_{12}$ extension.  The Galois group of $f_{E2}(t,x)$ 
  over $\Q(t)$ is
  $M_{12}.2$.  The $.2$ corresponds to the double cover of the $t$-line
  given by $z^2 = 11 (t-1)$.  
  
    The equation
  $f_{E2}(t,x)=0$ gives the genus zero degree twenty-four cover $X_{E2}$ of the $t$-line known to 
  exist by \cite[Prop.\ 9.1a]{MM}.   The curve $X_{E2}$ is just another name for the curve $X_E$ discussed above.
  It does not have any points over
  $\R$ or over $\Q_2$.    The 
  function $x$ has degree $2$, and there is a second function $y$
  so that $X_{E2}=X_E$ is given by $y^2=-x^2+40 x-404$.  
    
 \section{Dessins and generators}
 \label{dessins}
 
\begin{figure}[htb]
\includegraphics[width=4in]{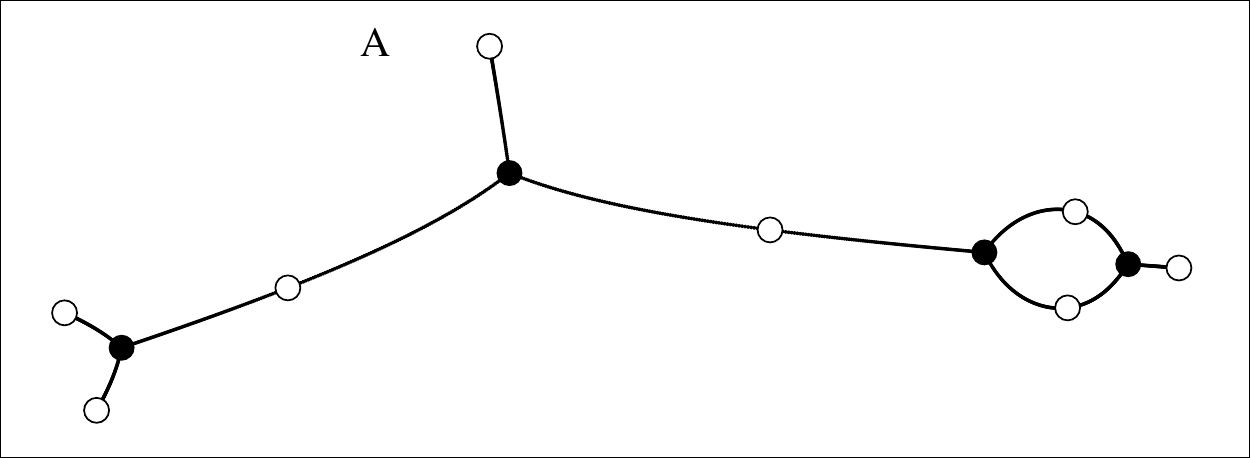}

\vspace{.2in}

\includegraphics[width=4in]{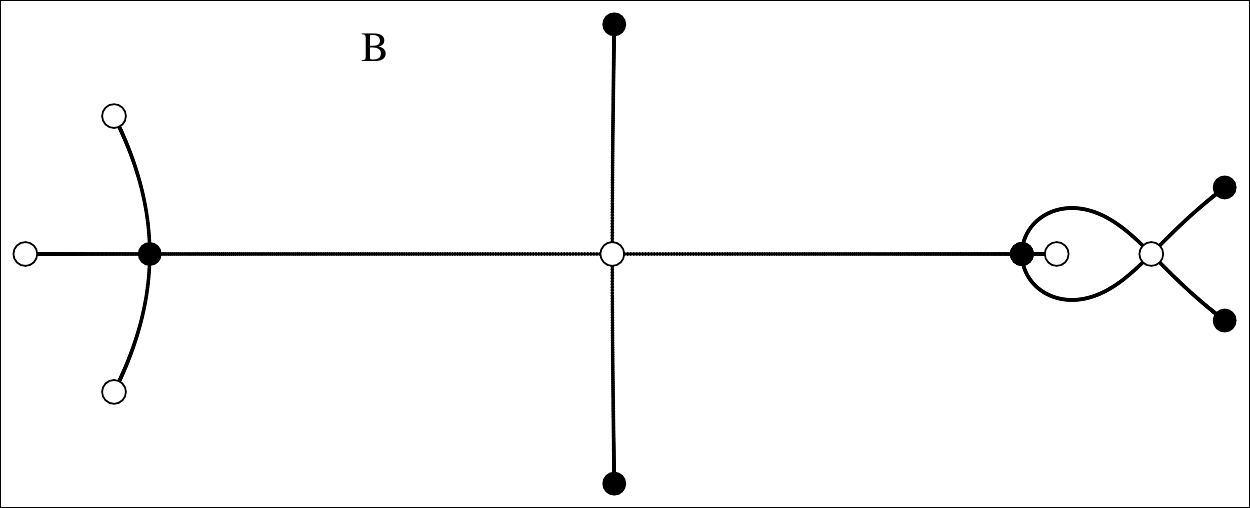}

\vspace{.2in}

\includegraphics[width=4in]{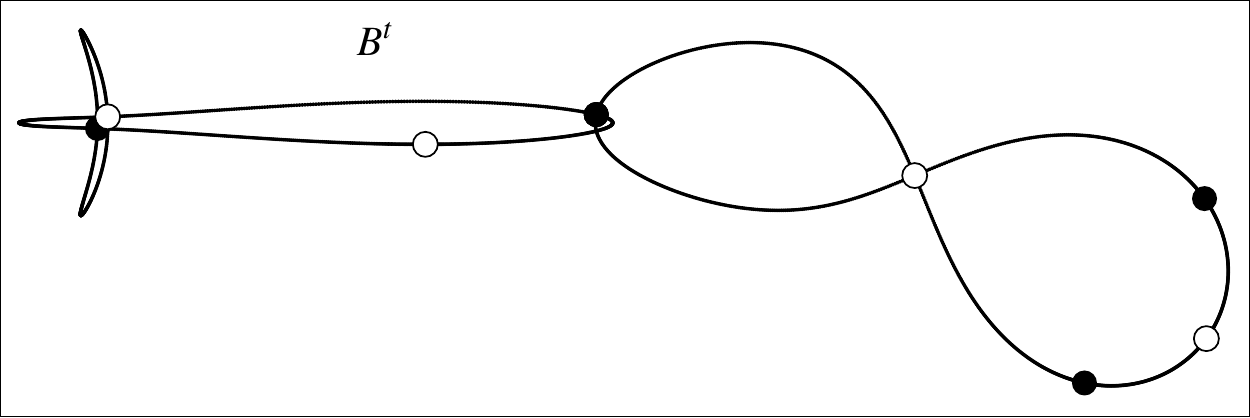}

\caption{\label{dessins5}  Dessins for Covers $A$, $B$, and $B^t$.  For Covers $A$ and $B$, the ambient
surface is the plane of the page; for Cover $B^t$ it is a genus two 
double cover of the plane of the page. }
\end{figure}

\begin{figure}[htb]
\[
\begin{array}{c}
\includegraphics[width=2.4in]{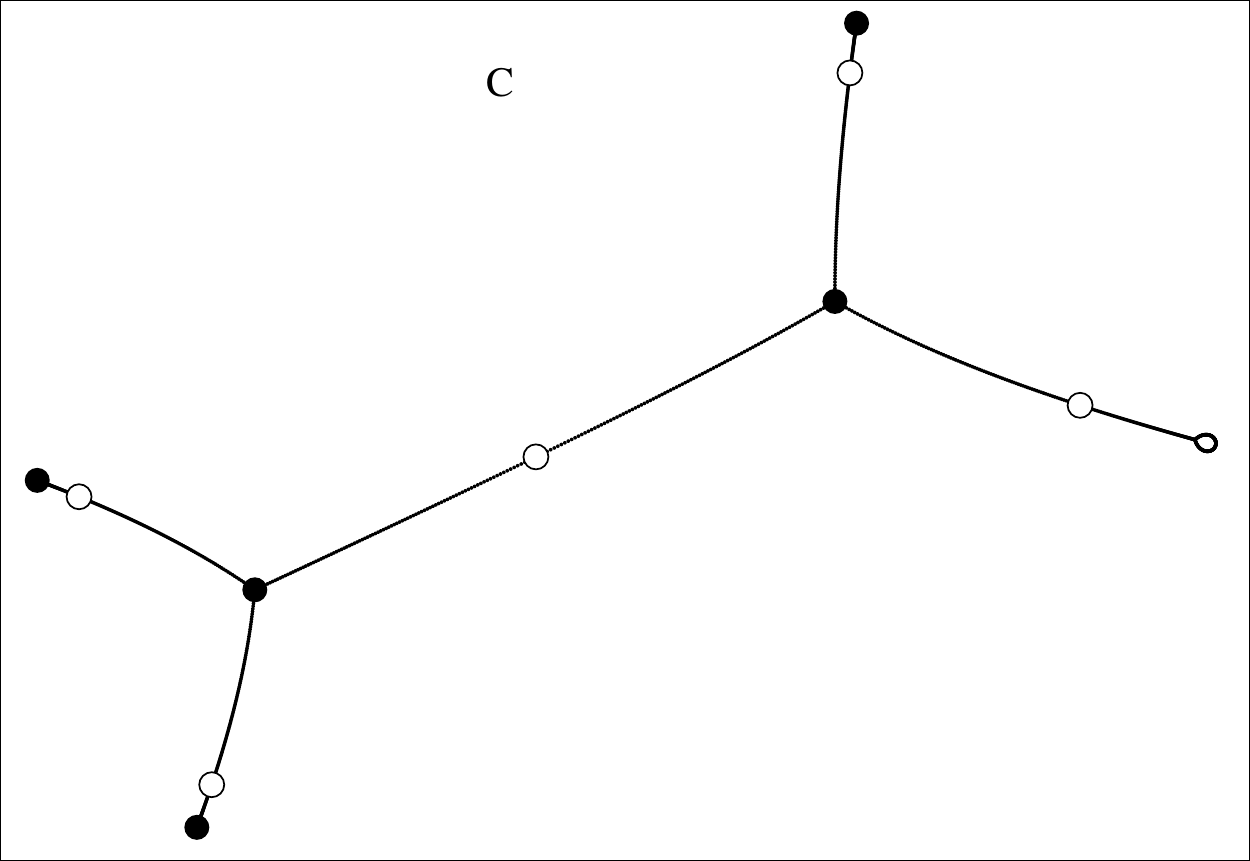}\\
  \\
\includegraphics[width=2.4in]{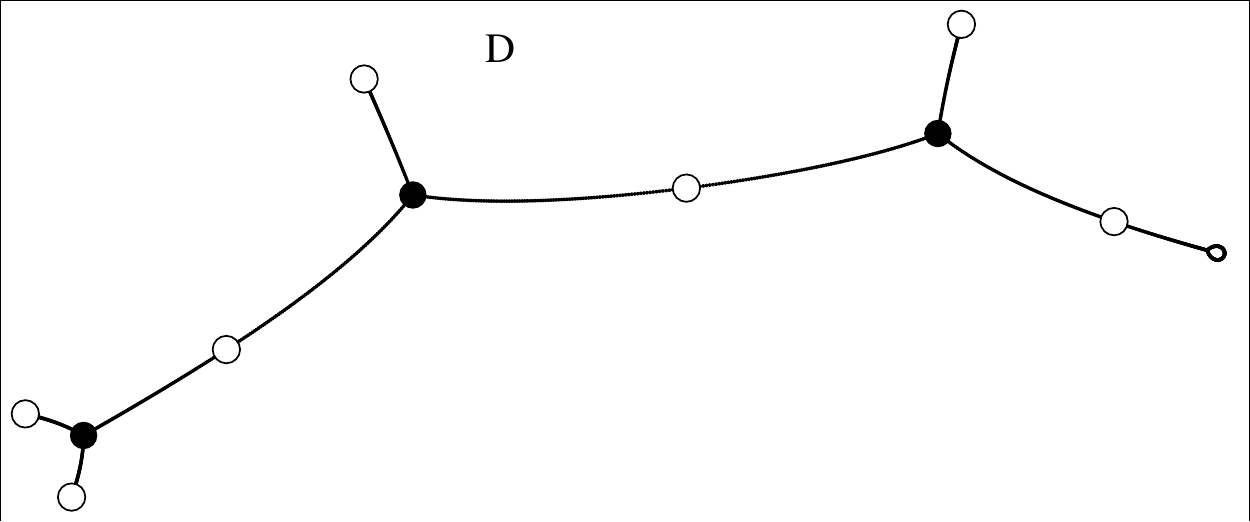} 
\end{array}
\;
\begin{array}{c}
\includegraphics[width=2.2in]{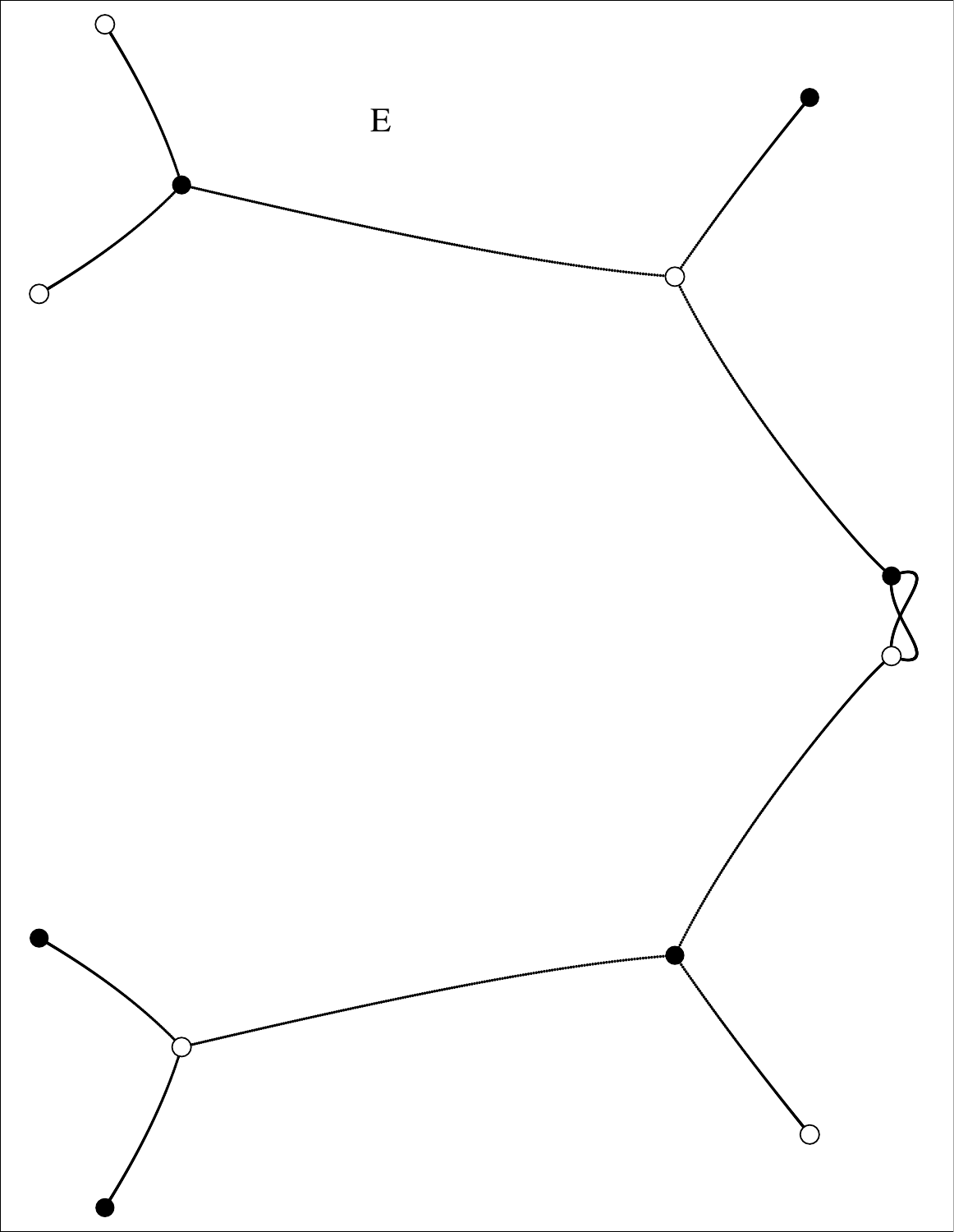}
\end{array}
\]

\caption{\label{dessins11} Dessins for Covers $C$, $D$, and $E$.   
For $C$ and $D$, the rightmost black and white vertices, of valence $3$ and $2$ respectively, are not drawn, so as not to obscure
the small loop to the right.  For Covers $C$ and $D$, the ambient surface is  the plane of the
page; for Cover $E$ it is a genus zero double cover of the plane of the page.   }
\end{figure}

     Figure~\ref{dessins5} draws pictures corresponding to Covers $A$, $B$,
 and $B^t$ while Figure~\ref{dessins11} draws pictures corresponding
 to Covers $C$, $D$, and $E$.   In this section, we describe the figures and
 how they give rise to generators of $M_{12}$ and $M_{12}.2$.   
 To avoid clutter, the twelve edges are not labeled in the figures.
 To follow the discussion, the reader needs to label the edges by $1$, \dots, $12$, in a way consistent with the text.

 \subsection{Covers $A$, $C$, and $D$}  For $L = A$, $C$, or $D$, 
 the corresponding figure draws the roots of $f_L(0,x)$ as black dots
 and the roots of $f_L(1,x)$ as white dots.   As $t$ moves from $0$ 
 to $1$, the roots of $f_L(t,x)$ sweep out the twelve edges of 
 the figure.   All together, the drawn bipartite graph, viewed as a subset
 of the Riemann sphere $X_L(\C)$, is the dessin of Cover $L$.  
 
     Let $T(\C)= \C-\{0,1\}$.  Let $\star = 1/2 \in T$.  
 Then the fundamental group $\pi_1(T(\C),\star)$ is 
 the free group $\langle \gamma_0,\gamma_1 \rangle$ where 
 $\gamma_k$ is the counterclockwise 
 circle of radius $1/2$ about $k$.   One
 has a natural extension $\pi_1(T(\C)_\R,\star)$ 
 obtained from $\pi_1(T(\C),\star)$ by adjoining
 a complex conjugation operator $\sigma$ 
 satisfying the involutory relation $\sigma^2 = 1$ and the
 intertwining relations $\sigma m_0 = m_{0}^{-1} \sigma$, 
 and $\sigma m_1 = m_{1}^{-1} \sigma$.
 
    For each $L$,  in conformity with the general theory of dessins, 
 we have a homomorphism 
$\rho_L$ from $\pi_1(T(\C),\star)$ into the group of permutations of the twelve edges.
Always the image is $M_{12}$.   Again in conformity
with the general theory, $\rho_L(\gamma_0)$ 
is the minimal counterclockwise rotation
about the black dots while $\rho_L(\gamma_1)$ 
is the minimal counterclockwise rotation
about the white dots.  

    Covers $A$, $C$, and $D$ behave very similarly.  
 Taking Cover $D$ as an example, and indexing  
 edges in the figure roughly from left to right, one has
 \begin{align*}
\rho_D(\gamma_0) & = (1,2,3)(4,5,6)(7,8,9)(10,11,12), \\
\rho_D(\gamma_1) & =  (3,4)(5,7)(8,10)(11,12).
\end{align*}
In the cases of $A$, $C$, and $D$, the twin permutation is easily visualized
as follows.   Consider the complex conjugates of the
dessins, obtained by 
flipping the drawn pictures upside down.  
In the new pictures, the image of edge $e$ is 
denoted $\overline{e}$.  Then we can apply
the general theory again.  In the case of Cover $D$
the result is  
\begin{align*}
 \rho^t_{{D}}(\gamma_0) & =
 (\overline{3},\overline{2},\overline{1}) 
 (\overline{6},\overline{5},\overline{4}) 
 (\overline{9},\overline{7},\overline{6}) 
 (\overline{12},\overline{11},\overline{10}), 
  \\
 \rho^t_{{D}}(\gamma_1) & =
 (\overline{3},\overline{4}) 
 (\overline{5},\overline{7}) 
 (\overline{8},\overline{10}) 
 (\overline{11},\overline{12}). 
\end{align*}
In other words, if one identifies the twelve $e$ with
their corresponding $\overline{e}$, 
one has $\rho_D(\gamma_0) = \rho_D^t(\gamma_0)^{-1}$ 
and $\rho_D(\gamma_1) = \rho_D^t(\gamma_1)$.  
The same relations hold with $D$ replaced by $A$ or $C$. 

The permutation representations $\rho_L$ and
$\rho^t_L$ are extraordinarily similar to each 
other in our three cases $L=A$, $C$, and $D$.  
Continuing with our example of Case $D$, 
considers words $w$ of length $\leq 14$ in 
$\gamma_0$ and $\gamma_1$.  Then 
$339$ different permutations arise 
as $\rho_D(w)$.  For each one of them,
the cycle type of $\rho_D(w)$ and
the cycle type of $\rho^t_D(w)$ agree.  
Only for words of length $15$ does one first get
a disagreement:
\begin{align*}
\rho_D(\gamma_0^2 \gamma_1 \gamma_0 \gamma_1 \gamma_0 \gamma_1 \gamma_0^2 \gamma_1 \gamma_0 \gamma_1 \gamma_0^2 \gamma_1) & = (1,8,6,5)(4,9,7,11)(2)(3)(10)(12), \\
\rho^t_D(\gamma_0^2 \gamma_1 \gamma_0 \gamma_1 \gamma_0 \gamma_1 \gamma_0^2 \gamma_1 \gamma_0 \gamma_1 \gamma_0^2 \gamma_1) & = (\overline{2},\overline{11},\overline{10},\overline{6})(\overline{3},\overline{4},\overline{8},\overline{9})(\overline{1},\overline{5})(\overline{7},\overline{12}).
\end{align*}
Only for words of length $18$ does one first get the other possible disagreement with 
$\rho_D(w)$ and $\rho_D^t(w)$ having different cycle types, one
$84$ and the other $821^2$.  

Reinterpreting, one immediately gets homomorphisms 
$\rho_{L2} : \pi_1(T(\C)_\R,\star) \rightarrow S_{24}$ with image
 $M_{12}.2$.  Namely the twenty-four element set is 
$\{1,\dots,12\} \cup \{\overline{1},\dots,\overline{12}\}$.  
The element $\rho_{L2}(\sigma)$ acts by 
interchanging each $e$ with its $\overline{e}$.   
For $k \in \{0,1\}$ one has $\rho_{L2}(\gamma_k) = \rho_{L}(\gamma_k) \rho_{L}^t (\gamma_k)$.

\subsection{Covers $B$ and $B^t$}  Inside the Riemann sphere $X_B(\C)$, we use black dots to represent
roots of $f_{B}(-\sqrt{5},x)$ and white dots to represent roots of $f_{B}(\sqrt{5},x)$.  
The twelve edges then correspond to the twelve preimages in the $x$-line 
of the interval $(-\sqrt{5},\sqrt{5})$.   Now we have monodromy
operators $m_\pm = m_{\pm \sqrt{5}}$ and a complex
conjugation $\sigma$ as before.  From the picture, indexing edges roughly from left to right again,
one immediately has
\begin{align*}
\rho_{B}(m_-) &= (1,2,4,3) (7,9,8,10) (5)(6)(11)(12), \\
\rho_{B}(m_+) & =  (4,5,7,6)(9,11,12,10)(1)(2)(3)(8), \\
\rho_{B}(\sigma) & =  (2,3)(5,6)(9,10)(11,12)(1)(4)(7)(8).
\end{align*}
The fact that Cover $B$ is defined over $\R$ corresponds to
$\rho_{B}(\sigma)$ already being in $M_{12}$.  

There are complications with presenting the twin case $B^t$
 visually, since $X_{B^t}$ has genus
two.  If we drew things using the $x$-variable, we would
have four black dots and four white dots,  all distinct 
in the $x$-plane.    An advantage of this presentation
would be that complex conjugation would be represented 
by the standard flip; this flip would fix exactly two 
of the black dots and all four of the white dots.  A 
disadvantage would be that some edges would be right
on top of other edges, and some edges would fold back
on themselves.  

Instead, we perturb things slightly, using the $y$-variable
of \S\ref{CovsB}, introducing the new variable $z = x + i y/200$,
and drawing the dessin instead in the $z$ plane.  Now monodromy
operators $\rho_{B^t}(m_-)$ and $\rho_{B^t}(m_+)$ can be
easily read off the picture.   Even $\rho_{B^t}(\sigma) = (2,3)(7,8)(9,10)(11,12)$
can be clearly read off,  the deformed version of the
real axis through all four white dots being easily imagined.

\subsection{Cover $E$}  Inside the sphere $X_E(\C)$, we draw the roots of $f_{E}(-\sqrt{-11},x)$ as 
black dots and the roots of $f_E(\sqrt{-11},x)$ as white dots.   As $s$ moves upward
from $-\sqrt{-11}$ to $\sqrt{-11}$, the twelve roots of $f_E(s,x)$ sweep out
the twelve drawn edges.   We work now with monodromy operators 
$m_\pm = m_{\pm \sqrt{-11}}$.    Because the ramification points 
$\pm \sqrt{-11}$ are no longer real, complex conjugation now satisfies
$\sigma m_- = m_+^{-1} \sigma$ and $\sigma m_+ = m_-^{-1} \sigma$.

The changes do not obstruct our basic procedure.  From the picture we have
\begin{align*}
\rho_E(m_+) & = \rho_E^t(m_-) = (3, 4, 5)(6, 7, 8)(10, 11, 12), \\
\rho_E(m_-) & = \rho_E^t(m_+) = (10,9,8)(5,6,7)(1,2,3).
\end{align*}
In this case, as just indicated, the twin representation is obtained by reversing the roles
of the black and white dots.  For the monodromy representation and its twin, 
 complex conjugation acts by subtraction from $13$ on
indices: $\rho_E(\sigma) = \rho^t_E(\sigma) = (1,12)(2,11)(3,10)(4,9)(5,8)(6,7)$.

The degree 24 dessin corresponding to $f_{E2}(t,x)$ can be easily imagined
from our drawn degree 12 dessin corresponding to $f_E(s,x)$.  Namely,
one first views both white dots and black dots as associated to the
number $t=0$.     Then one adds say a cross
at the appropriate midpoint $\epsilon$ of each edge $e$, viewing
these twelve crosses as associated to the number $t=1$.    
Any old edge $e$ now splits into two edges $\epsilon b$ and $\epsilon w$, with
$\epsilon b$ incident on a black vertex and $\epsilon w$ incident on a
white vertex.   The monodromy operator about zero is then
\begin{align*}
\rho_{E2}(m_0) & \! \! = \! \! (3b, 4b, 5b)(6b, 7b, 8b)(10b, 11b, 12b)(10w,9w,8w)(5w,6w,7w)(1w,2w,3w).
\end{align*}
The operator $\rho_{E2}(m_1)$ acts by switching each $\epsilon b$ and $\epsilon w$ while
the complex conjugation operator $\rho_{E2}(\sigma)$ acts by switching $\epsilon c$ and $(13-\epsilon)c$ for
 either color $c$.

\section{Specialization}
\label{specialization}

    This section still focuses on covers, but begins the process
 of passing from covers to number fields.  The next two sections are
 also focused on specialization, but with the emphasis shifted to the
 number fields produced.

\subsection{Keeping ramification within $\{2,3,q\}$} 
    Let $f(t,x) \in \Z[t,x]$ define a cover ramified only above the points $0$, $1$, and $\infty$ on
the $t$-line.  Then for each $\tau \in T(\Q) = \Q-\{0,1\}$, one has an associated number 
algebra $K_\tau$.  When $f(\tau,x)$ is separable, which it is for all our covers and all $\tau$, this
number algebra is simply $K_\tau = \Q[x]/f(\tau,x)$.   Thus ``specialization'' in our context
refers essentially to plugging in the constant $\tau$ for the variable $t$.  

\subsubsection*{Local behavior}
To analyze ramification in $\Q[x]/f(\tau,x)$, one works prime-by-prime.  
 The  procedure is described methodically in \cite[\S3,4]{ABC} and we review
it  in briefer and more informal language here.   For a given
 prime $p$, one puts $T(\Q)$ in the larger set $T(\Q_p) = \Q_p-\{0,1\}$.  One
 thinks of $T(\Q_p)$ as consisting of a generic ``center'' and three ``arms," one
 extending to each of the cusps $0$, $1$, and $\infty$.    A point 
 $\tau$ is in arm $k \in \{0,1,\infty\}$ if $\tau$ reduces to $k$ modulo $p$.
Otherwise, $\tau$ is generic.  If $\tau$ is in arm $k$, then one has its extremality
index $j \in \Z_{\geq 1}$, defined by $j = \ord_p(\tau)$,  $j=\ord_p(\tau-1)$, and
$j = - \ord_p(\tau)$ for $k=0$, $1$, and $\infty$ respectively.

    Suppose a prime $p$ is not in the bad reduction set of $f(t,x)$.   Then
 the analysis of $p$-adic ramification in any $K_\tau$ is very simple.
 First, if $\tau$ is generic, then $p$ is unramified in $K_\tau$.
 Second, suppose $\tau$ is in arm $k$ with extremality index $j$; then
 the $p$-inertial subgroup of the Galois group of $K_\tau$ is
 conjugate to $g^j_k$, where $g_k \in C_k$ is the
 local monodromy transformation about the cusp $k$.  In
 particular, suppose $g_k$ has order $m_k$;  then
 a point $\tau$ on the arm $k$ yields a $K_\tau$ 
 unramified at $p$ if and only if its extremality index
 is a multiple of $m_k$.  
    
  \subsubsection*{Global specialization sets}
    Let $S$ be the set of bad reduction primes of $f(t,x)$, thus $\{2,3,q\}$ for
us with $q=5$ for Covers $A$, $B$, and $B^t$ and $q=11$ for Covers $C$, $D$, and 
$E$.   Then the subset of $T(\Q)$ consisting of $\tau$ giving $K_\tau$ ramified only within $S$ depends only on $S$ and the monodromy orders $m_0$, $m_1$, and $m_\infty$.  
Following \cite{ABC} still, we denote it $T_{m_0,m_1,m_\infty}(\Z^S)$.
This set can be simply described without reference to $p$-adic numbers as follows.
It consists of all $\tau = - a x^p/c z^r$ where $(a,b,c,x,y,z)$ are 
integers satisfying the $ABC$ equation $a x^p+b y^q + c z^r=0$ 
with $a$, $b$, and $c$ divisible only by primes in $S$.  
After suitable simple normalization conditions are imposed, 
the integers $a$, $b$, $c$, $x$, $y$,
and $z$ are all completely determined by $\tau$.

    To find elements in some $T_{m_0,m_1,m_\infty}(\Z^S)$ one can 
carry out a computer search, restricting to $|ax^p|$ and $|cz^r|$ 
less than a certain height cutoff, say of the form $10^u$.  As one increases
$u$, the new $\tau$ found rapidly become more sparse.  
Many of the new $\tau$ are not entirely new, as they are 
often base-changes of lower-height $\tau$ as described in
\cite[\S4]{ABC}.    A typical situation, present for us here,
is that one can be confident that one has found at least most of
$T_{m_0,m_1,m_\infty}(\Z^S)$ by a short implementation
of this process.  

\begin{table}[htb]
\[
{
{\renewcommand{\arraycolsep}{3.4pt}
\begin{array}{rrclllllllll}
\mbox{$A2$:} &         |T_{3,2,10}^5| & = & 447,  && 158470321^3 & \!\! - \!\! & 1994904202391^2 
& \!\! + \!\! &  2^{10}  3^4 5^1 19^{10} & \!\! = \!\! & 0 \\
\mbox{$B$,$B^t$:} &\!\!\!  |T_{(4,4),10}^5| & = & 27, && 79^4 & \!\!-\!\! &  6881^2 &\!\! + \!\! &2^{8} 3^8 5 &\!\! = \!\! & 0 \\
\mbox{$C2$,$D2$:} & |T_{3,2,11}^{11}| & = & 394, & \!\!  \!\! & 2540833^3 &  \!\! - \!\!  &4050085583^2   &  \!\! + \!\!  &  2^{18} 3^1 11^6     & \!\!=\!\! & 0  \\
\mbox{$E2$:} & |T_{3,2,12}^{11}| & = & 395, && 796531585^3 & \!\!-\!\! & 22481204531903^2 & \!\!+\!\! &
 2^{11} 3^5 11^2 17^{12} & \!\!=\!\! & 0 \\
\end{array}
}
}
\]
\caption{\label{globspec} Sizes and largest height elements of specialization sets}
\end{table}
    The sizes of our specialization sets 
     $T_{m_0,m_1,m_\infty}^q \subseteq T_{m_0,m_1,m_\infty}(\Z^{\{2,3,q\}})$ are given by the left columns of Table~\ref{globspec}.  The right columns give the ABC triple corresponding to the element $\tau$ of largest  height in these 
sets.
The set $T^{5}_{(4,4),10}$ is not in our standard form.  We obtain it by considering
a set $T^5_{4,2,10}$ of $237$ points.  We select from this set the $\tau$ for which
$5(1-\tau)$ is a perfect square.   Each of these gives two specialization points 
$\sigma = \pm \sqrt{5(1-\tau)}$ in $T_{(4,4),10}$ and then 
we consider $\sigma = 0$ as in $T^{5}_{(4,4),10}$ as well.   The displayed
ABC triple yields $\sigma = \pm 6881/2^4 3^4$.

\begin{table}[htb]
{\small
{\renewcommand{\arraycolsep}{1.1pt}
\[
 \begin{array}{|c|c|cccccccccc|}
 \hline
 \mbox{gen} & \tau &    1 & 2 & 3 & 4 & 5 & 6 & 7 & 8 & 9 & 10  \\
     \hline
2&0&    (68) & &\multicolumn{3}{c}{\mbox{\bf Cover $A2$}} &&&&& \\  
  &1&  62 & (50) &  &  &  &  &  &  &  &    \\
   &\infty& 72 & \{46,52\} & 66 & \{46,48\} & 64 & 42 & 60 & \{40,42\} & 52 & 42  \\
      &  & 54  & (36 & 52  & 40 &  52 & 40  & 48  & 40 & 52  & 40 \\
      & &  52) & &  &  &  &  &  &  &  &   \\
      \hline
3&0&    52 & \{40,48\} & (36 & 48 & 48)&  &  &  &  &   \\
 & 1&    \{40,44\} & 36 & 24 & (20) &  &  &  &  &  &   \\
&\infty&    52 & 48 & 24 & 42 & 38 & 22 & 32 & 34 & 22 & 30 \\
& & 24 & 22 & 22 & 22 & (0 & 20 & 16 & 20 & 16 & 12   \\
    \hline
 5&0&      34 & 26 & (20 & 12 & 20) &  &  &  &  &    \\
&1&    34 & 26 & (18) &  &  &  &  &  &  &   \\
26,18 &\infty&    (42 & 42 & 42 & 42 & 26) &  &  &  &  &   \\
\hline
\multicolumn{12}{c}{\;} \\
\hline
2 &  0 & (\{18,20,24\}) & & \multicolumn{3}{c}{\mbox{\bf Cover $B$}} &&&&&\\
 &  1 &  (34)& &&&&&&&&\\
  & \infty  & \{16,22\} & 30 & \{16,22\} & 30 & \{12,18\} & 28 & \{12,16\} & 24 & \{12,18\} & 24 \\
   && (\{0,12\} & 22 & \{8,16\} & 22 & \{8,16\} & 18 & \{8,16\} &  22 & \{8,16\} & 22) \\
 \hline
 3 & \infty & 16 & 16 & 10 & 14 & 12 & 10 & 10 &( 10 & 0 & 10 \\
8,10 && 8 & 10 & 8 & 6 & 8 & 10 & 8) &  &  &   \\
 \hline
5     & 0 & 14 & (8) & &  & & & & & & \\
 18 & \infty & (10 & \{6,10\} & 20 & 18 & 20 & 18 & \{8,12\} & 18 & 20 & 18 ) \\
 \hline
   \end{array}
 \]
 }
 }
\caption{\label{abtable} Specialization tables for Covers $A2$ and $B$.}
\end{table}

\begin{table}[htb]
{\small
 \[
 {\renewcommand{\arraycolsep}{4pt}
 \begin{array}{|c|c|cccccccccc|}
 \hline
 \mbox{gen} & \tau &    1 & 2 & 3 & 4 & 5 & 6 & 7 & 8 & 9 & 10  \\
     \hline
   && && \multicolumn{3}{c}{\mbox{\bf Cover $C2$}} &&&&& \\
 2&0&  48 & \{12,24\} & 36 & 24 & 24 & (0 & 20 & 20) &  &    \\
  &1&  36 & (36 & 48) &  &  &  &  &  &  &    \\
   &\infty& 48 & \{12,24\} & 36 & 24 & 24 & 0 & 20 & 20 & 20 & 20  \\
   &        & 20     & 20    & 20  &  20 &  20 & 20)   &   &   &   &   \\
   \hline
   3&0& \{32,36\} & (36 & \{20,24\} & 36) & & &  &  & &
     \\
  &1&  34 & 22 & (20 & 16) &  &  &  &  &  &    \\
  24&\infty&  42 & 38 & \{18,22\} & 32 & 34 & 22 & 30 & 24 & (22 & 22 \\
   &        &   22   & 0    &  22 & 22  &22   &   22 & 22  &22   &22)   &   \\
      \hline
  11 & 0 &
    36 & 28 & (20 & 20 & 16) &  &  &  &  &   \\
   &1&  32 & (22) &  &  &  & &  &  &  &    \\
  24&\infty&  (44 & 44 & 44 & 44 & 44 & 44 & 44 & 44 & 44 & 44 \\
 &        &   24)   &     &   &   &   &    &   &   &   &   \\
 \hline
 \multicolumn{12}{c}{\;} \\
 \hline
  && && \multicolumn{3}{c}{\mbox{\bf Cover $D2$}} &&&&& \\
2 & 0 &   40 & \{16,24\} & 24 & 24 & 24 & (0 & 20 & 20) &  &      \\
  &1&  24 & (24 & 32) &  &  &  &  &  &  &    \\
  &\infty&  40 & \{16,24\} & 24 & 24 & 24 & (0 & 20 & 20 & 20 & 20  \\
   &        &  20    & 20    &  20 & 20   & 20  &  20)  &   &   &   &   \\
 \hline
3&0&    52 & \{40,48\} &(36 & 48 & 48) &  &  &  &  &    \\
  &1&  \{40,44\} & 36 & 24 & (20) &  &  &  &  &  &    \\
36,20&\infty&    52 & 48 & 24 & 42 & 38 & 22 & 32 & 34 & 22 & 30 \\
 &        &   24   & 22     &  22 &   22&  (0  & 20   & 20   & 20  & 20 & 20   \\
 &        &  20 & 20 & 20 & 20 & 20) & & &  & & \\
 \hline
11&0&    36 & 28 & (22 & 22 & 18) &  &  & &  &   \\
   &1& 32 & (20) &  &  &  &  &  &  &  &    \\
  24&\infty&   (44 & 44 & 44 & 44 & 44 & 44 & 44 & 44 & 44 & 44  \\
   &        &   24)    &     &   &   &   &    &   &   &   &   \\
 \hline
 \multicolumn{12}{c}{\;} \\
 \hline
   2 &0& 66 & 40 & 52 & \{24,32\} & 36 & 32 & 32 & (16 & 24 & 24) \\
  &1& 66 & (\{44,48\})&  \multicolumn{3}{c}{\mbox{\bf Cover $E2$}}  &  &  &  &  &    \\
     & \infty & 70 & (a & 74 & b & 72 & b & 74 & a &74 & b \\ 
 &        &  72    &  b   & 74)  &   &   &    &   &   &   &   \\
   \hline
 3     &0&48 & \{32,40\} & 32 & (40 & 40 & 32) &  &  &  &    \\
 &1&   48 & \{32,36\} & 32 & (24 & 28) &  &  &  &  &  \\
24,20&\infty&    56 & 52 & 32 & 48 & 48 &( 24,8 & 46 & 44 & 30 & 40 \\
 &        &  46    &  28   & 46  & 40  & 30   &  44  & 46)   &   &   &   \\
\hline
 11&0&  36 & 28 & (20 & 20 & 16) &  &  &  &  &   \\
  &1&    32 & 24 & (16 & 20) &  &  &  &  &  &   \\
  36,24 & \infty&  40 & 40 & 36 & 36 & 32 & 32 & 28 & 28 & 24 & 24 \\
 &        &   (0   & 22     & 20  & 18  & 16  &   22 &   12& 22   &  16 & 18  \\
 && 20 & 22) & &&&&&&&\\
\hline
   \end{array}
   }
 \]
 }
 \caption{\label{cdetable} Specialization tables for Covers $C2$, $D2$ and $E2$, 
with $a=\{24,36,48\}$ and $b = \{32,40,52\}$ in the case of Cover $E2$.}
 \end{table}
 
 \subsection{Analyzing $2$-, $3$-, and $q$-adic ramification}  Let $p \in \{2,3,q\}$.  
Then the quantity $\ord_p(\mbox{disc}(K_\tau))$ is a locally constant function on $T(\Q_p)$.  It 
shares some basic features with the much simpler tame case of 
$\ord_p(\mbox{disc}(K_\tau))$ for $p \not \in \{2,3,q\}$. For example, it is ultimately
periodic near each of the cusps.  However there are no strong general
theorems to apply in this situation, and the current best way to proceed is
computationally.

Each entry on Tables~\ref{abtable} and \ref{cdetable} gives a 
value of $\ord_p(K_\tau)$ for the indicated
cover and for $\tau$ in the indicated region.  The entries in the far left
column correspond to the generic region.  The entries in the main part of the
table correspond to the regions of the arms.  

For example, consider Cover $A2$ for $p=2$ and focus on the $\infty$-arm.  
This case is relatively complicated, as the table has three lines giving
entries corresponding to extremalities $1$-$10$ on the first line, $11$-$20$
on the second, and $21$ on the third.    A sample entry is $\{46,52\}$,
corresponding to extremality $j=2$.  This means first of all that
$\ord_2(K_\tau)$ can be both $46$ and $52$ in this region.   It
means moreover that our computations strongly suggest  that no other values
of $\ord_2(K_\tau)$ can occur.  The parentheses indicate the 
experimentally-determined periodicity.  Thus from the
table, $\ord_2(K_\tau) = 36$ is the only possibility for extremality  $12$, and
it is likewise the only possibility for extremalities $12 + 11 k$.  
We have no need of rigorously confirming the correctness of these 
tables, as they serve only as a guide for us in our search
for lightly ramified number fields.  Rigorous confirmations 
would involve computations which can be highly detailed
for some regions. Examples of interesting such computations
are in \cite{PV}.  

\subsection{Field equivalence.}  
\label{fe}A typical situation is that $f(\tau,x)$ is ireducible but 
has large coefficients.  
Starting in the next subsection, we apply {\em Pari}'s command {\em polredabs} \cite{Par} or some other procedure to obtain
a polynomial $\phi(x)$ with smaller coefficients defining the same field.  In general 
we say that two polynomials $f$ and $\phi$ in $\Q[x]$ are {\em field equivalent}, and write 
$f \approx \phi$, if $\Q[x]/f(x)$ and $\Q[x]/\phi(x)$ are isomorphic.   

\subsection{Specialization points with a Galois group drop}
\label{groupdrop}
 We now shift to explicitly indicating the source cover in the notation, writing $K(L,\tau)$ rather than $K_\tau$, as we will be
 often be considering various covers at once.  
Only a few of our algebras $K(L,\tau)$ have Galois group different
 from $M_{12}$ or $M_{12}.2$.   We present these degenerate  
 cases here,  before moving on to our main topic of non-degenerate
 specialization in the next section.

\begin{table}[htb]
 \[
 {\renewcommand{\arraycolsep}{2.1pt}
 {\small
 \begin{array}{|c|c|cc|ccc|ccc|rr|}
 \hline
  & &  && \multicolumn{3}{c|}{\mbox{Basic $p$-adic invariants}} &
   \multicolumn{3}{c|}{\mbox{Slope Content}}  & & \\
 \mbox{Cover} & \tau & \mbox{Fact} & G  & 2 & 3 & 5   & 2 & 3 & 5  & \mbox{RD} & \mbox{GRD}  \\
 \hline
\!\!\! B,B^t \!\!\! & \multirow{1}{*}{-5/2} &12 & L_2(11) &6_{8}^2  & 11_{10} 1 & 10_{13} 2 & 
[2]_3^2  &[\;]_{11}^5   & [\frac{3}{2}]_2& 41.2 & 55.4  \\
\hline
B &  \multirow{2}{*}{1}          & 11 \; 1 & M_{11}  &  6_{10} 4_8 1 \; 1 & 8_7 2_1\; 1 \; 1 & 5_9 \; 5_9 \;1 \; 1 & \multirow{2}{*}{ $[\frac{8}{3},\frac{8}{3}]_3^2$ } &  \multirow{2}{*}{$[\;]_8^2$} &  \multirow{2}{*}{$[\frac{9}{4}]_4$} & 96.2 &  \multirow{2}{*}{270.8}  \\
B^t & &   12 & M_{11}^t & 6_{10} \, 4_8 \, 2 & 8_7 \ 4_3 & 10_{19} \, 2_1    &  &&&103.3&   \\
\hline
B & \multirow{2}{*}{-11/5}  & 12 & M_{11}^t  & 6_{10} \, 4_8 \, 2  & 11_{10} 1 & 10_{19} \, 2  & 
\multirow{2}{*}{$[\frac{8}{3},\frac{8}{3}]_3^2$} &  \multirow{2}{*}{ $[\;]_{11}^5$   } &  \multirow{2}{*}{$[\frac{9}{4}]_4$} & 103.3 &  \multirow{2}{*}{281.2}  \\
B^t &  & 11 \; 1 & M_{11} &  6_{10} \, 4_8 \, 1 \; 1 & 11_{10} \, 1 & 9_5 \; 9_5 \; 1 \, 1 &  &&&117.5& \\
\hline
 \end{array}
 }
 }
 \]
 \caption{\label{groupdropb}  Description of $K_\tau$ and $K_\tau^t$ for the three $\tau$ in 
 $X_{(4,4),10}^5$ for which the Galois group is smaller than $M_{12}$} 
 \end{table}

For $B$ and $B^t$, our specialization set $T_{(4,4),10}^5$ has
twenty-seven points.  Three of them yield a group drop as in Table~\ref{groupdropb}.  
In this table, and also Tables~\ref{groupdropacde}, \ref{lowgrdm12}, \ref{lowgrdm122}, 
we present an analysis of ramification using the notation of \cite{JRLF} and
making use of the associated website repeatedly in the calculations.   
A $p$-adic field with degree $n = fe$, residual degree $f$, ramification
index $e$, and discriminant $p^{fc}$ is presented as $e^f_c$.  Superscripts $f=1$ 
are omitted.  Likewise subscripts $c = e-1$, corresponding to tame ramification, are
omitted.   Slope contents, as in  $[2]_3^2$, $[]_{11}^5$, and $[3/2]_2$ on the
first line, indicate decomposition groups and their natural filtration.  
This first field is tame at $3$ with inertia group of size $11$ and thus a contribution of $3^{10/11}$ to
the GRD.  It is wild at $2$ and $5$ with inertia groups of sizes $6$ and $10$ and contributions $2^{4/3}$ and
$5^{13/10}$ to the GRD respectively.  The Galois root discriminant, as printed, is $2^{4/3} 3^{10/11} 5^{13/10} \approx 55.4$.  

Continuing to discuss Table~\ref{groupdropb},
the specialization point $\tau = -5/2$ yields the same field in both $B$ and
 $B^t$, with group $L_2(11) = PSL_2(\F_{11})$ of
  order $660$.   A defining equation is 
  \begin{eqnarray*}
f_{B}(-5/2,x) & \approx & x^{12}-2 x^{11}-9 x^{10}+60 x^8+42 x^7+141 x^6+162 x^5+150 x^4 
\\ && \qquad +60
    x^3+141 x^2+18 x+21.
 \end{eqnarray*}
 The field $K(B,-5/2)$ is very lightly ramified, comparable with the remarkable
dodecic $L_2(11)$  field on \cite{KM} with  $\mbox{GRD} = \mbox{RD} = \sqrt{1831} \approx 42.8$.  
For the specialization point $\tau=1$, Cover 
$B$ yields a polynomial factorizing as $11+1$ while $B^t$ yields an irreducible
polynomial.   For the point $\tau = -11/5$ the situation is reversed.  Again
these fields are among the very lightest ramified of known fields with their
Galois groups, the first having been highlighted in our Section~\ref{background}.

For covers $A2$,  $C2$, $D2$, and $E2$ there are all together $1630$ specialization
points $\tau$.  Three of them yield group drops as in Table~\ref{groupdropacde}. 
\begin{table}[htb]
 \[
 {\renewcommand{\arraycolsep}{2.1pt}
 {\small
 \begin{array}{|c|c|cc|ccc|ccc|rr|}
 \hline
  & &  && \multicolumn{3}{c|}{\mbox{Basic $p$-adic invariants}} &
   \multicolumn{3}{c|}{\mbox{Slope Content}}  & & \\
 \mbox{Cover} & \tau & \mbox{Fact} & G  & 2 & 3 & 11   & 2 & 3 & 11  & \mbox{RD} & \mbox{GRD} \\
 \hline
C2  &-\frac{239^3}{3^{13}}  &24 & G_t & 2_3^6 2_3^6& 11_{10} 11_{10} 1\, 1 & 12_{11}^2 & 
 \multirow{2}{*}{$[2]_3^2$} &  \multirow{2}{*}{$[\;]_{11}^{10}$} & 
 \multirow{2}{*}{$[\;]_{12}^2$} & 63.6&  \multirow{2}{*}{$87.1$}  \\
   &  &12 & L_2(11).2 &2^6_3 & 11_{10}1 &12_{11} & 
&   &&63.6  &  \\
\hline
C2 &    \frac{3 \cd 11^5}{2^7}  & 22 \; 2 & G_i  &11_{10}^2  & 6_{11}6_{10}3_53_52_12 &4_3^24_34_32_1^22_1 & \multirow{2}{*}{$[\;]_{11}^{10}$ } &  \multirow{2}{*}{$[\frac{5}{2}]_2^2$} & 
 \multirow{2}{*}{$[\;]_4^2$} & 47.6&  \multirow{2}{*}{$85.0$}  \\
 & &   12 & L_2(11).2 & 11_{10} 1  & 6_{11}^2   &  4_3^2 4_3   &  &&&80.7 &   \\
\hline
D2 &  \frac{-17^3}{2^7} & 22 \; 2& G_i & 11_{10}^2   & 6_7 6_6 3_3 3_3 2_1 2 & 10_9 10_9 2_1    & 
\multirow{2}{*}{$[\;]_{11}^{10}$} &  \multirow{2}{*}{$[\frac{3}{2}]_2^2$} &  \multirow{2}{*}{$[\;]_{10}$} &
 40.4  &  \multirow{2}{*}{58.6}  \\
  &  & 12 & L_2(11).2 &  11_{10} 1 & 6_7^2 & 10_9 1 \, 1  &  &&&38.8& \\
\hline
 \end{array}
 }
 }
 \]
 \caption{\label{groupdropacde} Description of $K(L,\tau)$ for the only three 
 instances  where the  Galois group is smaller than $M_{12}$ in Cases
 $A2$, $C2$, $D2$, and $E2$  }
 \end{table}
 In all three cases, the Galois group is $PGL_2(11) = L_2(11)$, in either a transitive or an intransitive
 degree twenty-four representation.  The least ramified case is the last one, for which
 a degree twelve polynomial is 
\begin{equation*}
\mbox{{\small $
x^{12}-6 x^{10}-6 x^9-6 x^8+126 x^7+104 x^6-468 x^5+258 x^4
 +456 x^3-1062 x^2+774
    x-380.$}
    }
\end{equation*}
The GRD here is small, but still substantially larger than the smallest known GRD for a $PGL_2(11)$ number field
of $3^{10/11} 227^{1/2} \approx 40.90$.  This field comes from a modular form of weight one and conductor $3 \cdot 227$ in
characteristic $11$ \cite[App.~A]{Sch}.   The examples of this section serve to calibrate expectations for
the proximity to minima of the $M_{12}$ and $M_{12}.2$ number fields in the next section.

\section{Lightly ramified $M_{12}$ and $M_{12}.2$ number fields}
\label{light}      
   This section reports on ramification of specializations to fields ramified
   within $\{2,3,q\}$ with $q = 5$ for covers $A2$, $B$, $B^t$ and 
   $q=11$ for covers $C2$, $D2$, $E2$.    Our presentation
   continues to use the  conventions of \S\ref{fe} on field equivalence and of \S\ref{groupdrop} 
   on $p$-adic ramification.

   According to Tables~\ref{abtable} and \ref{cdetable}
 the maximal root discriminants our covers 
can give for these fields are
\begin{align*}
\delta_{A2}^{\rm max}  = 
(2^{72} 3^{52} 5^{42})^{1/24} & \approx  1445, &   \delta_{C2}^{\rm max} = (2^{48} 
3^{42} 11^{44})^{1/24} & \approx 2219, \\
\delta_{B}^{\rm max}  = (2^{34} 3^{16} 5^{18})^{1/12} & \approx \;\;  344, & \delta_{D2}^{\rm max} =  (2^{40} 3^{52} 11^{44})^{1/24} & \approx 2784, \\
&& \delta_{E2}^{\rm max} = (2^{74} 3^{56} 11^{40})^{1/24}  & \approx 5985.
\end{align*}
The fields highlighted below all have substantially smaller root discriminant.  
Subsections \S\ref{rd}, \S\ref{grd}, and \S\ref{ps} focus respectively on
fields with small root discriminant, small Galois root discriminant, and 
at most two ramifying primes.  

\subsection{Small root discriminant}
\label{rd}  The smallest root discriminant appearing for our $M_{12}$ specializations is approximately $46.2$, as reported on the first line
of Table~\ref{lowgrdm12} below.  This is substantially above the smallest known root discriminant $2^2 3^1 29^{1/3} \approx 36.9$ from
\cite{KM}, discussed above in \S\ref{pursuing}.  For the larger group $M_{12}.2$, the two smallest root discriminants appearing
in our list 
 are
 $(2^{12} 3^{24} 11^{22})^{1/24} \approx 38.2$ and 
$(2^{20} 3^{24} 11^{20})^{1/24} \approx 39.4$.  The smallest root discriminant 
comes from Cover $C2$ at $\tau = 5^3/2^2$ and the field can be given
by the polynomial
\begin{eqnarray*}
\lefteqn{f_{C2}(5^3/2^2,x) \approx} \\
&& x^{24}-11 x^{23}+53 x^{22}-154 x^{21}+330 x^{20}-594 x^{19}+1012
    x^{18}-2255 x^{17} \\ && +6512 x^{16} -17710 x^{15}+42768 x^{14}-89067
    x^{13}+154308 x^{12}-237699 x^{11} \\ && +351252 x^{10} -483318 x^9+623997
    x^8-753291 x^7+733491 x^6-520641 x^5 \\ &&+278586 x^4-104841 x^3 +15552 x^2+2673
    x+81.
 \end{eqnarray*}
 The second smallest root discriminant also comes from Cover $C2$.
 It arises twice, once from $-17^3/2^7$ and once from $7^3/2^9$. Both
 these specialization points define the same field.   There are 
 seven more $M_{12}.2$ fields with root discriminant under
 $50$, each arising exactly once.  In order, they come from
 the covers $D2$, $A2$, $D2$, $C2$, $A2$, $A2$, and $D2$.

\subsection{Small Galois root discriminant}
\label{grd}
 For $M_{12}$, the smallest known Galois root discriminant appears in \cite{KM} and also on the first
line of Table~\ref{lowgrdm12}.   
The fact that $E$ appears only once in Table~\ref{lowgrdm12} is just a reflection 
of the  simple fact that $q=5$ for $B$ and $B^t$ while $q=11$ for $E$.  
\begin{table}[htb]
{\small
 {\renewcommand{\arraycolsep}{3pt}
 \[
 \begin{array}{|c|c|ccc|ccc|rr|}
 \hline
  & & \multicolumn{3}{c|}{\mbox{Basic $p$-adic invariants}} &
   \multicolumn{3}{c|}{\mbox{Slope Content}}  & & \\
 \mbox{Cover} & \tau & 2 & 3 & q&  2 & 3 & q  & \mbox{RD} & \mbox{GRD}  \\
 \hline
 B & \multirow{2}{*}{5} & 8_{16} 3_2 1 & 11_{10} & 10_{13} 2_1 &
  \multirow{2}{*}{$[\frac{4}{3},\frac{4}{3},3]_3^2$} &  \multirow{2}{*}{$[\;]_{11}^5$} &
    \multirow{2}{*}{$[\frac{3}{2}]_2$} & 46.2 & \multirow{2}{*}{93.2} \\ 
B^t &                                  &8_{16} 4_4 & 11_{10} & 10_{13} 2_1 & &&&51.6 & \\
   \hline
    B & \multirow{2}{*}{0} & 12_{34} & 9_9 2_1 1& 3_2^23_2^2 &
  \multirow{2}{*}{$[\frac{23}{6},\frac{23}{6},3,\frac{8}{3},\frac{8}{3}]_3^2$} &  \multirow{2}{*}{$[\frac{9}{8},\frac{9}{8}]_{8}^2$} &  \multirow{2}{*}{$[\;]_3^2$} & 52.1 & \multirow{2}{*}{112.0} \\ 
B^t &                                  &12_{34} & 12_{12} & 3_2^23_2^2 & &&&62.5 & \\
   \hline
 B & \multirow{2}{*}{$5/2$} & 12_{12} & 9_9 2_1 1 & 10_{13} 2_1 & 
 \multirow{2}{*}{$[\frac{8}{3},\frac{8}{3},\frac{4}{3},\frac{4}{3}]_3^2$} &
  \multirow{2}{*}{$[\frac{9}{8},\frac{9}{8}]_8^2$} &
   \multirow{2}{*}{$[\frac{3}{2}]_2$} &
   58.2 &  \multirow{2}{*}{$132.4$} \\
 B^t &  & 12_{12} & 12_{12} &  10_{13} 2_1 & & & & 69.9 &   \\
 \hline
 B & \multirow{2}{*}{$-5$} & 4_8^3 & 9_9 2_1 1 & 10_{13} 2_1 & 
  \multirow{2}{*}{$[3,\frac{5}{2},2,2]^6$} &
  \multirow{2}{*}{$[\frac{9}{8},\frac{9}{8}]_8^2$} &
   \multirow{2}{*}{$[\frac{3}{2}]_2$} &
   65.3 &  \multirow{2}{*}{$153.0$} \\
 B^t &                                  & 4_8^3 & 12_{12} & 10_{13} 2_1 & & & & 78.5 &  \\
 \hline 
  B & \multirow{2}{*}{$-3$} & 8_{16} 4_4 & 8_7 2_1 1 \, 1 & 5_9^2 1^2 & 
  \multirow{2}{*}{$[3,\frac{4}{3},\frac{4}{3}]_3^2$} &
  \multirow{2}{*}{$[\;]_8^2$} &
   \multirow{2}{*}{$[\frac{9}{4}]_4^2$} &
   73.8 &  \multirow{2}{*}{$255.6$} \\
 B^t &                                  & 8_{16} 3_2 1 & 8_7 4_3& 10_{19} 2_1& & & & 103.3 &  \\
 \hline
 E & -319/54& 2_2^3 2_2^3 & 3_3^3 3_3 & 11_{16} 1 &
 \multirow{2}{*}{$[2]^3$} & 
  \multirow{2}{*}{$[\frac{3}{2}]^3$} & 
   \multirow{2}{*}{$[\frac{8}{5}]_5$} & 
   146.8 & \multirow{2}{*}{280.6} \\
 E &319/54& 2_2^3 2_2^3 & 3_3^3 3_3 & 11_{16} 1 & &&& 146.8 & \\
 \hline 
   B & \multirow{2}{*}{$-5/3$} & 6_{10} 4_8 2_2 & 9_{16} 1^2 1 & 10_{13} 2_1 & 
  \multirow{2}{*}{$[\frac{8}{3},\frac{8}{3},2]_3^2$} &
  \multirow{2}{*}{$[2,2]^2$} &
   \multirow{2}{*}{$[\frac{3}{2}]_2$} &
  89.8 &  \multirow{2}{*}{$287.9$} \\
 B^t &                                  & 6_{10} 4_8 2_2 & 9_{16} 1^2 1 & 10_{13} 2_1 & & & & 89.8 &  \\
 \hline 
 \end{array}
 \]
 }
 }
 \caption{\label{lowgrdm12} The fourteen $M_{12}$ fields from our list with Galois root discriminant
 $\leq 300$, grouped in twin pairs.  The two $\tau$'s for $E$ both come from 
 $\sigma = 23^3/2^2 3^6$. }
 \end{table}

For $M_{12}.2$, Galois root discriminants can be substantially smaller than the 
minimum known for $M_{12}$, 
as illustrated by Table~\ref{lowgrdm122}.  
In this case, in contrast to $M_{12}$, the field giving the smallest known root discriminant also
gives the smallest known Galois root discriminant. 
\begin{table}[htb]
{\small
 {\renewcommand{\arraycolsep}{0.7pt}
 \[
 \begin{array}{|c|c|ccc|ccc|rr|}
 \hline
  & & \multicolumn{3}{c|}{\mbox{Basic $p$-adic invariants}} &
   \multicolumn{3}{c|}{\mbox{Slope Content}}  & & \\
 \mbox{Cover} & \tau & 2 & 3 & q&  2 & 3 & q  & \mbox{RD} & \mbox{GRD}  \\
 \hline
 C2 & 5^3/2^2 & 3_2^6 1^6 & 9_{12} 3_3^2 3_3^2 1 \, 1 \,  1 & 12_{11}^2 &
 [\;]_3^6 & [\frac{3}{2},\frac{3}{2}]_2^2  & [\;]_{12}^2 &  38.2  & 65.8 \\
 C2 & 11^3/2^3 & 4_6^6 & 10_0 10_9 2_1 2_1 & 12_{11} 6_5 4_3 2_1 & [2,2]^4 & [\;]_{10}^4 & [\;]_{12}^2 & 52.1 & 68.5  \\ 
 C2 & -11^2/2^6 3^3 &  &12_{12} 9_9 2_1 1 & 22_{27} 2_1 &&   [\frac{9}{8},\frac{9}{8}]_8^2 & [\frac{13}{10}]_{10} & 63.3 & 69.1 \\ 
 A2 & -2^3 5^4 11^3/3^8 & 8_{18}  4^4_8 & 12_{18} 9_{15} 2_1 1 & 4_2 4_2 4_2 4_2 4_ 2 4_2 & 
 [3,2,2]^4 & [\frac{15}{8},\frac{15}{8}]_8^2 & [\;]_2^2 & 44.9 & 73.9 \\  
 C2 & \left\{ \begin{array}{c} -17^3/2^7 \\ 7^3/2^9 \end{array}  \right.& 11_{10}^2 1^2 & 9_{12} 3_3^2 3_3^2 1 \, 1\, 1 & 11_{10} 11_{10} 2_1 2_1 & [\;]_{11}^2 & [\frac{3}{2},\frac{3}{2}]_2^2 & [\;]_{10} & 39.4 & 74.7 \\ 
 A2 & -5^4/2^33^3 & 8_{22}^2 8_{22} & 9_{12} 9_9 3_3 3 & 4_24_24_24_24_24_2 & [\frac{7}{2},3,2,2]^2 &
  [\frac{3}{2},\frac{3}{2}]_2^3 & [\;]_2^2 &  45.1 &75.4 \\
C2 &  5^3/3^3 &  2_3^6 2_3^6 & 12_{12} 9_9 2_1 1 & 10_910_92_12_1 & [3]^6 & [\frac{9}{8},\frac{9}{8}]_8^2 & [\;]_2 & 57.1 & 81.7 \\
 C2 & -2^9 5^3/3^2 & & 9_{18} 6_{10} 6_{10} 1^2 & 12_{11}^2 & & [\frac{9}{4},\frac{9}{4}]_4^2 & [\;]_{12}^2 & 51.3 & 88.8 \\
 D2 & 5^3/2^2 & 3_2^4 3_2^4 & 9_{15} 6_9 3_4 3_3 2_1 1 & 12_{11} 6_5 4_3 2_1 & [\;]_3^4 & [2,\frac{3}{2}]_2^2 & [\;]_{12}^2 & 50.7 & 94.8 \\
 D2 & -11^2/2^6 3^3 && 12_{12} 9_{12} 1^2 1 & 22_{27} 2_1 && [\frac{3}{2},\frac{3}{2}]_2^4 & [\frac{13}{10}]_{10} & 49.2  & 95.2 \\
 \hline
 \end{array}
 \]
 }
 }
 \caption{\label{lowgrdm122} The ten $M_{12}.2$ fields from our list with the Galois root discriminant $\leq 100$.}
 \end{table}

\subsection{At most two ramifying primes}
\label{ps}

Let $d_L(\tau)$ be the field discriminant of $K(L,\tau)$.  Then, generically for our specializations,
$d_L(\tau)$ has the form $\pm 2^a 3^b q^c$ with all three exponents positive.  The few cases 
where at least one of the exponents is zero are as follows.
For Cover~$A2$, from Table~\ref{abtable} the prime $3$ drops out from the discriminant 
exactly if $\ord_3(\tau) \in \{-15, -25, -35, \dots\}$.  This drop occurs in $2$ of our $447$ specializations:
\begin{eqnarray*}
d_{A2}(71^3/2^3 3^{15} 5^2) & = & 2^{66} 5^{42},  \\
d_{A2}(3289^3/2^7 3^{15} 5) & = & 2^{60} 5^{42}. 
\end{eqnarray*}
For Covers~$C2$ and $D2$, from Table~\ref{cdetable} the prime $2$ drops out exactly if $\ord_2(\tau) \in \{6,9,12,\dots\} \cup \{6,17,28,\dots\}$.  
This much less demanding condition is met by $24$ of our $394$ specialization points, as listed in Table~\ref{2drop}.
\begin{table}[htb]
\[
   \begin{array}{rclrr}
   &\!\!\!\!\!\! \tau \!\!\!\!\!\! & & d_{C2}(\tau) & d_{D2}(\tau) \\
   \hline
    -2^9 5^3  & \fr & 3^2 & 3^{38} 11^{22} & 3^{48}
      11^{20} \\
    -11^2 & \fr & 2^6 3^3 & 3^{22} 11^{28} & 3^{24}
      11^{28} \\
    -2^9 & \fr & 3^3 & 3^{22} 11^{32} & 3^{24} 11^{32}
      \\
    -11 \cd 131^3 & \fr & 2^6 3^3 & 3^{18} 11^{36} &
      3^{24} 11^{36} \\
    -3^3 11 & \fr & 2^6 & 3^{20} 11^{36} & 3^{36}
      11^{36} \\
    -11 & \fr & 2^6 & 3^{34} 11^{36} & 3^{44}
      11^{36} \\
    2^6 11 & \fr & 3^6 & 3^{22} 11^{36} & 3^{22}
      11^{36} \\
    11 \cd 59^3  & \fr & 2^{17} 3^2 & 3^{38} 11^{36} &
      3^{48} 11^{36} \\
    -67^3 & \fr & 2^6 11 & 3^{34} 11^{44} & 3^{40}
      11^{44} \\
    -2^{12} & \fr & 11 & 3^{34} 11^{44} & 3^{40}
      11^{44} \\
    -2^6 3^3 & \fr & 11^2 & 3^{20} 11^{44} & 3^{36}
      11^{44} \\
    -2^6 & \fr & 11 & 3^{34} 11^{44} & 3^{40}
      11^{44} \\
       \end{array}
      \;\;\;\;\;\;\;\;\;\;
    \begin{array}{rclrr}
   &\!\!\! \!\!\!\tau \!\!\! \!\!\!& & d_{C2}(\tau) & d_{D2}(\tau) \\
   \hline
    -2^6 31^3  & \fr & 3^3 11^5 & 3^{22} 11^{44} &
      3^{20} 11^{44} \\
    -2^6 & \fr & 3^3 11 & 3^{18} 11^{44} & 3^{24}
      11^{44} \\
    173^3 & \fr & 2^6 11^7 & 3^{34} 11^{44} & 3^{40}
      11^{44} \\
    2^9 & \fr & 3 \cd 11^3 & 3^{42} 11^{44} & 3^{52}
      11^{44} \\
    13^3 & \fr & 2^6 11^2 & 3^{34} 11^{44} & 3^{44}
      11^{44} \\
    7^3 & \fr & 2^6 11 & 3^{24} 11^{44} & 3^{20}
      11^{44} \\
    2087^3 & \fr & 2^6 3^{15} 11 & 3^{20} 11^{44} &
      11^{44} \\
    3^6 & \fr & 2^6 11 & 3^{24} 11^{44} & 3^{28}
      11^{44} \\
    553^3 & \fr & 2^6 3^9 11^2 & 3^{22} 11^{44} &
      3^{22} 11^{44} \\
    313^3 & \fr & 2^6 3^6 11 & 3^{22} 11^{44} & 3^{22}
      11^{44} \\
    89^3 & \fr & 2^6 11^2 & 3^{24} 11^{44} & 3^{32}
      11^{44} \\
    7033^3 & \fr & 2^6 3^6 11^4 & 3^{22} 11^{44} &
      3^{22} 11^{44}
   \end{array}
\] 
\caption{\label{2drop} The $24$ specialization points of $T_{3,2,11}^{11}$ at which
$2$ does out from discriminants of specializations in the $C2$ and $D2$ families}
\end{table}

Also for Covers $C2$ and $D2$ it is possible for $3$ to drop out.  From Table~\ref{cdetable} this occurs if
$\ord_3(\tau)$ is in $\{-12,-23,-34,\dots\}$ for $C2$ or in $\{-15,-26,-37,\dots\}$ for $D2$.  
This stringent condition is met once in each case.  For $C2$, this one
$3$-drop gives $d_{C2}(-55177^3/2^3 3^{23} 11^2) = 2^{36} 11^{44}$.    
For $D$, the one $3$-drop occurs where there is
also a $2$-drop.  An equation for the corresponding number field 
is given at the end of \S\ref{lifts24}.

\section{Lifts to the double covers $\tilde{M}_{12}$ and $\tilde{M}_{12}.2$}
\label{lifts}

    In this final section, we discuss lifts to $\til{M}_{12}$ and $\til{M}_{12}.2$.  
Interestingly, our six cases behave quite differently from each other.   

\subsection{Lack of lifts to $(2.M_{12}.2)^*$}
\label{lack}
The $.2$ for the geometrically disconnected degree twenty-four covers $A2$, $C2$, and $D2$
corresponds to the constant imaginary quadratic fields $\Q(\sqrt{-5})$, 
$\Q(\sqrt{-11})$, and $\Q(\sqrt{-11})$ respectively.   Accordingly
all specializations have complex conjugation in the class $2C$ on
Table~\ref{conjclasses}.  Elements of $2C$ lift to 
elements of order $4$ in the nonstandard $(2.M_{12}.2)^*$ as reviewed
in \S\ref{isocline}.  Thus $M_{12}.2$ fields of the form $K(L2,\tau)$ 
with $L \in \{A,C,D\}$ do not embed in $(2.M_{12}.2)^*$ fields.  

The $B$ families do not even give $M_{12}.2$ fields.
Also $M_{12}.2$ fields of the form $K(E2,\tau)$ do not
embed in $(2.M_{12}.2)^*$ fields, as explained at the end of 
\S\ref{lifts12} below.   For these reasons, we have 
deemphasized $(2.M_{12}.2)^*$ in this paper, 
despite the fact that it fits into the framework of our title.
An open problem which we do not pursue here is 
to explicitly write down a degree forty-eight polynomial
in $\Q[y]$ with Galois group $(2.M_{12}.2)^*$.   

\subsection{Geometric lifts to $\til{M}_{12}$}   In this subsection, we work over $\C$ so that symbols
such as $X_L$ should be understood as complex algebraic curves.  Table~\ref{lifttable} reprints the 
six partition triples belonging to $M_{12}$ of Table~\ref{sixtriples} and for each indicates
lifts to partition triples in $\til{M}_{12}$.    For a fixed label $L$ in 
$\{A,B,B^t,C,D,E\}$, let 
$(g_0,g_1,g_\infty)$ be the permutation triple in $M_{12}$ presented pictorially
in Figure~\ref{dessins5} or \ref{dessins11}.  
 Then our main focus is a permutation triple $(\tilde{g}_0,\tilde{g}_1,\tilde{g}_\infty)$.
Here each $\tilde{g}_k$ is in the class indicated in row $\til{L}$ and column $k$ of
Table~\ref{lifttable}.    
One has $\tilde{g}_0 \tilde{g}_1 \tilde{g}_\infty = 1$ and accordingly one
gets a double cover $\til{X}_{{L}}$ of $X_L$.  The degree $24$ map 
$\til{X}_{{L}} \rightarrow \bbP^1$ by design has monodromy group
$\til{M}_{12}$.  

The class of $-\tilde{g}_k$ is indicated on Table~\ref{lifttable} right below the
class of $\tilde{g}_k$.   Note that $\pm \tilde{g}_0$, $\pm \tilde{g}_1$, and
$\pm \tilde{g}_\infty$ multiply either to $1$ or $-1$ in $\til{M}_{12}$,
according to whether the number of minus signs is even or odd.   Thus
our choice of $(\tilde{g}_0,\tilde{g}_1,\tilde{g}_\infty)$ could equally
well be replaced by $(\tilde{g}_0,-\tilde{g}_1,-\tilde{g}_\infty)$,
$(-\tilde{g}_0,\tilde{g}_1,-\tilde{g}_\infty)$, or 
$(-\tilde{g}_0,-\tilde{g}_1,\tilde{g}_\infty)$.  
The choice we make always minimizes the
genus of $\tilde{X}_L$.   

\begin{table}[htb]
\[
\begin{array}{|ccccc|c|ccccc|}
\cline{1-5} \cline{7-11}
\mbox{Cover} & 0 & 1 & \infty & g && \mbox{Cover} & 0 & 1 & \infty & g \\
\cline{1-5} \cline{7-11}
\multicolumn{11}{c}{\;} \\
\cline{1-5} \cline{7-11}
A & 3^4 & 2^4 1^4 & (10) 2 & 0    && C & 3^3 1^3 & 2^6 & (11)1 & 0 \\
\til{A}    & 3^8 & 2^8 1^8 & (20) 4 & 0   && \til{C}     &3^61^6 & 4^6 & 11^2 1^2 & 2  \\
    & 6^4 & 2^{12} & (20)4 &&&&  6^3 2^3 & 4^6 & (22)2  &   \\
\cline{1-5} \cline{7-11}   
\multicolumn{11}{c}{\;} \\
\cline{1-5} \cline{7-11}
B & 4^21^4 & 4^2 1^4 & (10) 2 & 0 && D & 3^4 & 2^4 1^4 & (11) 1 & 0 \\
\til{B} &  4^4 2^2 1^4 & 4^4 2^2 1^4 & (20) 4  & 2 && \til{D} & 3^8 & 2^81^8 & (22)2 & 0 \\
                  & 4^4 2^2 1^4 & 4^4 2^2 1^4 & (20) 4 &     &&                 & 6^4 & 2^{12} &  11^2 1^2 &  \\
\cline{1-5} \cline{7-11}
\multicolumn{11}{c}{\;} \\
\cline{1-5} \cline{7-11}
B^t & 4^2 2^2 & 4^2 2^2 & (10) 2 & 2 && E & 3^3 1^3 & 3^3 1^3 & 66 & 0 \\
\til{B}^t & 4^4 2^4 & 4^42^4 & (20)4 & 4 & & \til{E} & 3^6 1^6 & 3^6 1^6 & 12^2 & 0 \\
                 & 4^4 2^4 & 4^42^4 & (20)4 &     & &               & 6^3 2^3 & 6^3 2^3 & 12^2 & \\
\cline{1-5} \cline{7-11}
\end{array}
\]
\caption{\label{lifttable} Lifts of partition triples in $M_{12}$ to partition triples
in $\til{M}_{12}$.}
\end{table}

\begin{figure}[htb]
\includegraphics[width=5in]{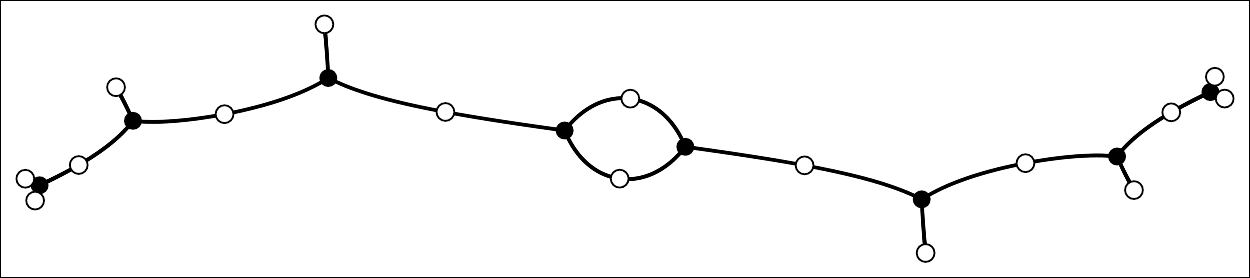}
\caption{\label{pictddisks2} The dessin of $f_{\tilde{D}}(t,y)$ in the complex $y$-line} 
\end{figure}

To understand Table~\ref{lifttable} in diagrammatic terms, consider Cover $D$ as an 
example.   The curve $X_D$ is just the complex $x$-line.  Its cover $\tilde{X}_D$
is just the complex $y$-line, with relation given by $y = x^2$.   The dessin
drawn in $X_D$ in Figure~\ref{dessins11} ``unwinds'' to a double
cover dessin in $\tilde{X}_D$ drawn in Figure~\ref{pictddisks2}.    The monodromy operators
are
\begin{align*}
\tilde{g}_0 = \rho_{\tilde{D}}(\gamma_0) & = (1,2,3)(4,5,6)(7,8,9)(10,11,12) && \\ &  \qquad (-1,-2,-3)(-4,-5,-6)(-7,-8,-9)(-10,-11,-12)  \\
\tilde{g}_1 = \rho_{\tilde{D}}(\gamma_1) & = (3,4)(5,7)(8,10)(11,-12)(-11,12)(-3,-4)(-5,-7)(-8,-10)
\end{align*}
Cover $A$ is extremely similar, with $\tilde{X}_A$ also covering $X_A$ via
$y=x^2$.  

At a geometric level, the six covers behave similarly, as just described.  
However
The curves $X_L$ are defined over $\Q(\sqrt{-5})$, 
$\Q$, $\Q$, $\Q(\sqrt{-11})$,
$\Q(\sqrt{-11})$, and $\Q$ for $L = A$, $B$, $B^t$, $C$, $D$, and $E$ 
respectively.  At issue is whether the cover $\til{X}_L$ can likewise
be defined over this number field.

\subsection{Lifts to $\tilde{M}_{12}$}
\label{lifts12}
{\em A lifting criterion.}  A $v$-adic field $K_v$ has 
a local root number $\epsilon(K_v) \in \{1,i,-1,-i\}$.  For example, 
taking $v=\infty$, one has $\epsilon(\C)=-i$; also if $K_p/\Q_p$ is unramified
then $\epsilon(K_p)=1$.   The invariant $\epsilon$
extends to algebras by multiplicativity: $\epsilon(K_v' \times K_v'') = \epsilon(K'_v)\epsilon(K_v'')$.
For an algebra $K_v$, there is a close relation between the local root
number $\epsilon(K_v)$ and the Hasse-Witt invariants $HW(K_v) \in \{-1,1\}$.  
In fact if the discriminant of $K_v$ is trivial as an element of $\Q_v^\times/\Q_v^{\times 2}$ then
$\epsilon(K_v) = HW(K_v)$.    In this case of trivial discriminant class,
one has $\epsilon(K_v)=1$ if and only if the homomorphism 
$h_v : \Gal(\overline{\Q}_v/\Q_v) \rightarrow A_n$ corresponding
to $K_v$ can be lifted into a homomorphism into the Schur double
cover, $\tilde{h}_v : \Gal(\overline{\Q}_v/\Q_v) \rightarrow \til{A}_n$.
If $K$ is now a degree $n$ number field then one has local root numbers 
$\epsilon(K_v)$ multiplying to $1$.  In the case when the discriminant
class is trivial, then all $\epsilon(K_v)$ are $1$  
if and only if the homomorphism 
$h : \Gal(\overline{\Q}/\Q) \rightarrow A_n$ corresponding
to $K$ can be lifted into a homomorphism
 $\tilde{h} : \Gal(\overline{\Q}/\Q) \rightarrow \til{A}_n$.   The general
 theory of local root numbers is presented in more detail
 in \cite[\S3.3]{JRLF} and local root numbers are calculated automatically
 on the associated database.  
 
 Since the map $\til{M}_{12} \rightarrow M_{12}$ is
 induced from $\til{A}_{12} \rightarrow A_{12}$,
 one gets that an $M_{12}$ number field embeds in
 an $\til{M}_{12}$ number field if and only
 if all $\epsilon(K_v)$ are trivial.   Also it 
 follows from the above theory that if $K$ and $K^t$ are twin
 $M_{12}$ fields then $\epsilon(K_v) = \epsilon(K_v^t)$ for 
 all $v$.  
 
 \subsubsection*{Covers $B$ and $B^t$.}
From the very last part of \cite{Me}, summarizing the 
approach of \cite{BLV}, we have the general formula
\begin{equation}
\label{hilbertsymbol}
\epsilon(K(B,\tau)_v) = (25-5 \tau^2,\tau)_v,
\end{equation}
where the right side is a local Hilbert symbol.  
For example, one gets
that $K(B,\tau)$ is obstructed at $v=\infty$
if and only if $\tau < -\sqrt{5}$.  Similar
explicit computations identify exactly 
the locus of obstruction for all primes $p$.
This locus is empty if and only if $p \equiv 3, 7 \; (20)$.  

In particular, because one has obstructions even in specializations,
$\tilde{X}_B$ cannot be defined over $\Q$.  One can see this
more naively as follows.    By Table~\ref{lifttable}, one 
has eight points on $X_B$ corresponding to the $1^4 1^4$.  
The cover $\tilde{X}_B$ is ramified at exactly four of these points,
those which correspond to the $2^2 2^2$.   But the eight points in 
$X_B$
are the roots of
\[
 x^8+36 x^7+462 x^6+2228 x^5-585 x^4-30948 x^3-22388 x^2+215964 x-82539
 \]
 and this polynomial is irreducible in $\Z[x]$.  
 
 If $\tau \in \Q$ is a square then of course all the local Hilbert symbols
 in \eqref{hilbertsymbol} vanish. This motivates consideration
 of the following base-change diagram of smooth projective
 complex algebraic curves:
 \[
 \begin{array}{ccccccl}
 \til{Z}_B & \rightarrow & \til{X}_B  && 7 && 2  \\
\downarrow & & \downarrow && && \\ 
 Z_B & \rightarrow & X_B &\;\;\;\;\;\; \mbox{(with genera} &0 && 0 \;\;\; ).  \\
 \downarrow &  & \downarrow && &&  \\
 \bbP^1 & \rightarrow & \bbP^1 &&0& &0 \\
 \end{array}
 \]
 Here the bottom map is the double cover $t \mapsto t^2$ 
 and $Z_B$ and $\tilde{Z}_B$ are the induced double 
 covers of $X_B$ and $\tilde{X}_B$ respectively.  Thus
 the cover $Z_B \rightarrow \bbP^1$ is a five-point cover,
 with ramification invariants $2^4 1^4$ above
 the fourth roots of $5$ and $5^2 1^2$ above
 $\infty$.    While $\til{X}_B$ 
 is not realizable over $\Q$, Mestre proved that 
 $\til{Z}_B$ is realizable \cite{Me}.  
 
 Rather than seek equations for a degree $24$ polynomial
 giving the cover $\tilde{Z}_B$, we content ourselves with
 a single example in the context of number fields.   Conveniently
 the first twin pair on Table~\ref{lowgrdm12} is unobstructed.  
 Corresponding equations are 
\begin{eqnarray*}
f_{B}(5,x) & \approx &  x^{12}-2 x^{11}+6 x^{10}+15 x^8-48 x^7+66 x^6-468 x^5-810 x^4 \\ && +900 x^3+486
    x^2+1188 x-1314, \\
    &&\\
f_{B^t}(5,x) & \approx & x^{12}-2 x^{11}+6 x^{10}+30 x^9-30 x^8+60 x^7-150 x^6+120 x^5-285 x^4 \\ && +150
    x^3-120 x^2+90 x+30.
\end{eqnarray*}
Carefully taking square roots of the correct field elements, to avoid Galois groups such as
the generically occuring
$2^{12}.M_{12}$, we find covering $\tilde{M}_{12}$ fields to be given by
\begin{eqnarray*}
\tilde{f}_{B}(5,y) & \approx & y^{24}-30 y^{20}+540 y^{18}+945 y^{16}-22500 y^{14}-58860 y^{12}+421200
    y^{10} \\ && +1350000 y^8-7970400 y^6+11638080 y^4-6480000 y^2+1166400, \\
    &&\\
\tilde{f}_{B^t}(5,y) & \approx & y^{24}+40 y^{22}+480 y^{20}-1380 y^{18}-46260 y^{16}-10800 y^{14} \\ &&+1190340
    y^{12} -4429800 y^{10}+65650500 y^8-324806400 y^6 \\ && 
   +588257280 y^4-398131200
    y^2 +58982400.
\end{eqnarray*}
The $p$-adic factorization partitions of these polynomials 
 for the first $|\til{M}_{12}| = 190080$ primes different from $2$, $3$, and $5$ 
 are summarized in Table~\ref{conjclasses}.  As expected  
from the Chebotarev density theorem, the distribution is quite similar to the distribution
of elements of $\til{M}_{12}$ in classes.  The one class not represented is the central
non-identity class $1A2$.  Calculating now with five times as many primes, exact equidistribution
would give five classes each for $1A1$ and $1A2$.  In fact, in this range 
there are eight primes splitting at the $M_{12}$ level, 
$\mathit{76493}$, $\mathit{2956199}$, $\mathit{5095927}$, $7900033$, $\mathit{7927511}$, $\mathit{10653197}$, $11258593$, and $\mathit{12420649}$.  Those in ordinary type
correspond to $1A2$ while those in italics to $1A1$.  

\subsubsection*{Cover $E$}  For all $\tau \in \Q$, the algebra $K(E,\tau)$ is obstructed
at $\infty$, since $K(E,\tau)_\infty \cong \C^6$ and 
$\epsilon(\C^6) = \epsilon(\C)^6 = (-i)^6 = -1$.   This obstruction 
can be seen more directly from Table~\ref{conjclasses}: a field 
in $K(E,\tau)$ has complex conjugation in class $2A$ of $M_{12}$ of
cycle type $2^6$.  The only class above $2A$ in $\til{M}_{12}$ 
is $2A2$ of cycle type $4^6$, and so the complex conjugation element
cannot lift.

\subsection{Lifts to $\tilde{M}_{12}.2$}  
\label{lifts24}  Lifting for the remaining cases behaves as follows.  

\subsubsection*{Cover $A$.}   The polynomial $f_A(t,x)$ from \S\ref{CoverA} gives an equation 
for $X_A$.  From Table~\ref{lifttable}, we see that $\tilde{f}_A(t,x) = f_A(t,y^2)$ gives an equation
for $\tilde{X}_A$ with coefficients in $\Q(\sqrt{-5})$.   The Galois group of $f_A(t,y^2)$
over $\Q(\sqrt{-5})(t)$ is $\tilde{M}_{12}$ by construction.  However the
Galois group of the rationalized polynomial $f_{A2}(t,y^2)$ over $\Q(t)$ is not $\tilde{M}_{12}.2 = 2.M_{12}.2$ 
but rather an overgroup of shape $2^2.M_{12}.2$, with the 
final $.2$ corresponding to $\Q(\sqrt{-5})$ present already
in the splitting field of $f_{A2}(t,x)$.  

The overgroup also has shape $2.M_{12}.2^2$.  Here the
quotient $2^2$ corresponds to $\Q(\sqrt{3},\sqrt{-15})$.  
Over $\Q(\sqrt{3})$, the polynomial $f_{A2}(t,x^2)$
has Galois group $2.M_{12}.2$.  Over $\Q(\sqrt{-15})$ it
has Galois group the isoclinic variant $(2.M_{12}.2)^*$ discussed
in \S\ref{isocline}.

\subsubsection*{Cover $C$.}  Here Table~\ref{lifttable} says that $\tilde{X}_C$ is a double
cover of $X_C$ ramified at six points and hence of genus two.   A defining polynomial
is 
 \[
 \tilde{f}_C(t,y)  =  \mbox{Resultant}_x(y^2-2 h(x),f_C(t,x)) 
 \]
 where
 \begin{eqnarray*}
 h(x) & = & 2 x^6+22 x^5 u-22 y^5-165 x^4 u-957 x^4-1804 x^3 z+4664 x^3 \\ 
 && +4884 x^2 u+17754 x^2+4686 x u-15114 x+385 u+1243.
 \end{eqnarray*}
 Here the Galois group of the rationalized polynomial $\tilde{f}_{C2}(t,y)$ is indeed the desired $\tilde{M}_{12}.2$.  As an 
 example of an interesting specialization, consider $\tau = 5^3/2^2$ corresponding
 to the first line of Table~\ref{lowgrdm122}.    A corresponding polynomial is 
 \begin{eqnarray*}
\lefteqn{ \tilde{f}_{C2}(5^3/2^2,y)  \approx} \\
&&  y^{48}-22 y^{44}+495 y^{40}-4774 y^{36}+51997 y^{32}-214038 y^{28}+64152
    y^{26} \\ &&+2194852 y^{24} -705672 y^{22}-4044304 y^{20}-30696732
    y^{18}+61713630 y^{16} \\ && 
 +149602464 y^{14}   -9212940 y^{12}+569477304
    y^{10}+138870369 y^8 \\ && -484796664 y^6 +1029399030 y^4
    +39870468 y^2+793881.
 \end{eqnarray*}
 The fields $K(C2,5^3/2^2)$ and $\til{K}(C2,5^3/2^2)$ respectively have discriminant,
 root discriminant, and Galois root discriminant as follows:
 \begin{align*}
 d & = 2^{12} 3^{24} 11^{22}, & \tilde{d} & = 11^2 d^2,  \\
 \delta & =  2^{1/2} 3^1 11^{11/12} \;\;\;\;\,\,\, \approx 38.2, & \tilde{\delta} & = 11^{1/24} \delta  \,\, \approx 42.2,\\
 \Delta & =  2^{2/3} 3^{25/18} 11^{11/12} \approx 65.8, & \tilde{\Delta} & = 11^{1/24} \Delta  \approx 72.7.  
 \end{align*}
The first two splitting primes for $\tilde{f}_{C2}(5^3/2^2,x)$ are $1270747$ and $2131991$.

 The specialization point $\tau = 5^3/2^2$ just treated is well-behaved as follows.  
In general, to keep ramification of $\tilde{K}(C,\tau)$ within $\{2,3,11\}$, one must take specialization points
in the subset $T_{3,4,11}(\Z^{\{2,3,11\}})$ of 
$T_{3,2,11}(\Z^{\{2,3,11\}})$.     While the known part of 
$T_{3,2,11}(\Z^{\{2,3,11\}})$  has $394$ points, the subset in $T_{3,4,11}(\Z^{\{2,3,11\}})$ has only
$78$ points.  In particular, while $\tau = 5^3/2^2 \in T_{3,4,11}(\Z^{\{2,3,11\}})$, the other 
six specialization points for Cover $C2$ appearing in Table~\ref{lowgrdm122}
are not.

\subsubsection*{Cover $D$}   From Table~\ref{lifttable} we see that $f_{D}(t,y^2)=0$ is
an equation for $\tilde{X}_D$ with $\Q(\sqrt{-11})$ coefficients.   This equation combines
the good features of the cases just treated.  Like $\tilde{X}_A$ but unlike $\tilde{X}_C$,
the cover $\tilde{X}_D$ has genus zero.  Like Case $C$ but unlike
Case $A$, the Galois group of the rationalized polynomial $f_{D2}(t,y^2)$
 over $\Q(t)$ is  $\tilde{M}_{12}.2$.

    At the 2-3-dropping specialization point $2087^3/2^63^{15}11$ of Table~\ref{2drop},
a defining polynomial with $e = 11$ is as follows:
\begin{eqnarray*}
\lefteqn{\tilde{f}_{D2}(2087^3/2^6 3^{15} 11,y) \approx} \\
&&     y^{48} + 2 e^3 y^{42} + 69 e^5 y^{36} + 868 e^7 y^{30} -4174 e^7
    y^{26} + 11287 e^9 y^{24}  \\ && -4174 e^{10} y^{20} + 5340 e^{12} y^{18} + 131481
    e^{12} y^{14} +17599 e^{14} y^{12} + 530098 e^{14} y^8 \\ &&  + 3910 e^{16}
    y^6 + 4355569 e^{14} y^4 + 20870 e^{16} y^2 + 729 e^{18}.
\end{eqnarray*}

\noindent
The $p$-adic factorization patterns for the first $|\tilde{M}_{12}.2| = 380160$ primes
different from $11$ are summarized in Table~\ref{conjclasses}.   Again one sees
agreement with the Haar measure on conjugacy classes.  In this case, the first 
primes split at the $M_{12}.2$ level are
$3903881$, $8453273$, $11291131$, $12153887$, $15061523$, 
$15359303$.  Two of these are still split at the $\tilde{M}_{12}.2$
level, namely $11291131$ and $15061523$.  

   The Kl\"uners-Malle database \cite{KM} contains an $M_{11}$ field ramified 
at $661$ only.  The polynomial just displayed makes $M_{12}$ the second
sporadic group known to appear as a subquotient of the Galois
group of a field ramified at one prime only.   These two examples are
quite different in nature, because $661$ is much too big to divide $|M_{11}|$ 
while $11$ divides $|M_{12}|$.

 \end{document}